\documentclass[fleqn,11pt]{article}
\usepackage{mycommands,amssymb,amsmath,amsthm,color,pagesize,outlines,cite,caption}
\usepackage[pdftex]{graphicx}
\usepackage[round]{natbib}
\usepackage{hyperref} 
\usepackage{setspace}
\DeclareMathOperator*{\diag}{diag}
\DeclareMathOperator*{\Tr}{Tr}

\begin{document}

\newtheorem{Theorem}{Theorem}[section]
\newtheorem{Lemma}[Theorem]{Lemma}
\newtheorem{Corollary}[Theorem]{Corollary}
\theoremstyle{definition} \newtheorem{Definition}[Theorem]{Definition}

\title{Robust estimation of principal components from depth-based multivariate rank covariance matrix}
\date{}
\author{Subhabrata Majumdar\\
University of Minnesota, Twin Cities\\
\vspace{1em}
Technical Report}
\maketitle

Abstract:
Analyzing principal components for multivariate data from its spatial sign covariance matrix (SCM) has been proposed as a computationally simple and robust alternative to normal PCA, but it suffers from poor efficiency properties and is actually inadmissible with respect to the maximum likelihood estimator. Here we use data depth-based spatial ranks in place of spatial signs to obtain the orthogonally equivariant Depth Covariance Matrix (DCM) and use its eigenvector estimates for PCA. We derive asymptotic properties of the sample DCM and influence functions of its eigenvectors. The shapes of these influence functions indicate robustness of estimated principal components, and good efficiency properties compared to the SCM. Finite sample simulation studies show that principal components of the sample DCM are robust with respect to deviations from normality, as well as are more efficient than the SCM and its affine equivariant version, Tyler's shape matrix. Through two real data examples, we also show the effectiveness of DCM-based PCA in analyzing high-dimensional data and outlier detection, and compare it with other methods of robust PCA.
\vspace{.5cm}


\newpage
\section{Introduction}
In multivariate analysis, the study of principal components is important since it provides a small number of uncorrelated variables from a potentially larger number of variables, so that these new components explain most of the underlying variability in the original data. In case of multivariate normal distribution, the sample covariance matrix provides the most asymptotically efficient estimates of eigenvectors/ principal components, but it is extremely sensitive to outliers as well as relaxations of the normality assumption. To address this issue, several robust estimators of the population covariance or correlation matrix have been proposed which can be used for Principal Components Analysis (PCA). They can be roughly put into these categories: robust, high breakdown point estimators that are computation-intensive \citep{rousseeuw85, maronna76}; M-estimators that are calculated by simple iterative algorithms but do not necessarily possess high breakdown point \citep{huber77, tyler87}; and symmetrised estomators that are highly efficient and robust to deviations from normality, but sensitive to outliers and computationally demanding \citep{dumbgen98, sirkia07}.

When principal components are of interest, one can also estimate the population eigenvectors by analyzing the spatial sign of a multivariate vector: the vector divided by its magnitude, instead of the original data. The covariance matrix of these sign vectors, namely Sign Covariance Matrix (SCM) has the same set of eigenvectors as the covariance matrix of the original population, thus the multivariate sign transformation yields computationally simple and high-breakdown estimates of principal components \citep{locantore99, visuri00}. Although the SCM is not affine equivariant, its orthogonal equivariance suffices for the purpose of PCA. However, the resulting estimates are not very efficient, and are in fact asymptotically inadmissible \citep{magyar14}, in the sense that there is an estimator (Tyler's M-estimate of scatter, to be precise) that has uniformly lower asymptotic risk than the SCM.

The nonparametric concept of data-depth had first been proposed by \cite{tukey75} when he introduced the halfspace depth. Given a dataset, the depth of a given point in the sample space measures how far inside the data cloud the point exists. An overview of statistical depth functions can be found in \citep{zuo00}. Depth-based methods have recently been popular for robust nonparametric classification \citep{jornsten04, ghosh05, dutta12, sguera14}. In parametric estimation, depth-weighted means \citep{ZuoCuiHe04} and covariance matrices \citep{ZuoCui05} provide high-breakdown point as well as efficient estimators, although they do involve choice of a suitable weight function and tuning parameters. In this paper we study the covariance matrix of the multivariate rank vector that is obtained from the data-depth of a point and its spatial sign, paying special attention to its eigenvectors. Specifically, we develop a robust version of principal components analysis for elliptically symmetric distributions based on the eigenvectors of this covariance matrix, and compare it with normal PCA and spherical PCA, i.e. PCA based on eigenvectors of the SCM.

The chapter is arranged in the following fashion. Section \ref{section:sec2} provides preliminary theoretical concepts required for developments in the subsequent sections. Section \ref{section:dcmSection} introduces the Depth Covariance Matrix (DCM) and states some basic results related to this. Section \ref{section:sec4} provides asymptotic results regarding the sample DCM, calculated using data depths with respect to the empirical distribution function, as well as its eigenvectors and eigenvalues. Section \ref{section:simSection} focuses solely on principal component estimation using the sample DCM. We obtain influence functions and asymptotic efficiencies for eigenvectors of the DCM. We also compare their finite sample efficiencies for several multinormal and multivariate $t$-distributions with those of the SCM, Tyler's scatter matrix and its depth-weighted version through a simulation study. Section \ref{section:sec6} presents two applications of the methods we develop on real data. Finally, we wrap up our discussion in Section \ref{section:sec7} by giving a summary of our findings and providing some potential future areas of research. Appendices \ref{section:appA} and \ref{section:appB} contain all technical details and proofs of the results we derive.

\section{Preliminaries}\label{section:sec2}
\subsection{Spatial signs and sign covariance matrix}
Given a vector $\bfx \in \mathbb{R}^p$, its spatial sign is defined as the vector valued function \citep{MottonenOja95}:
$$ \bfS(\bfx) = \begin{cases} \bfx\| \bfx \|^{-1} \quad \mbox{if }\bfx \neq \bf0\\
\bf0 \quad \mbox{if }\bfx = \bf0 \end{cases} $$
When $\bfx$ is a random vector that follows an elliptic distribution $|\Sigma|^{-1/2} f((\bfx - \bfmu)^T \Sigma^{-1} (\bfx - \bfmu))$, with a mean vector $\bfmu$ and covariance matrix $\Sigma$, the sign vectors $\bfS(\bfx - \bfmu)$ reside on the surface of a $p$-dimensional unit ball centered at $\bfmu$. Denote by $\Sigma_S(\bfX) = E\bfS (\bfX - \bfmu)\bfS (\bfX - \bfmu)^T$ the covariance matrix of spatial signs, or the \textit{Sign Covariance Matrix} (SCM). The transformation $\bfx \mapsto \bfS(\bfx - \bfmu)$ keeps eigenvectors of population covariance matrix unchanged, and eigenvectors of the sample SCM $ \hat \Sigma_S = \sum_{i=1}^n \bfS (\bfx_i - \bfmu)\bfS (\bfx_i - \bfmu)^T/n $ are $\sqrt n$-consistent estimators of their population counterparts \citep{taskinen12}.

The sign transformation is rotation equivariant, i.e. $ \bfS(P (\bfx - \bfmu)) = P(\bfx - \bfmu)/\| P (\bfx - \bfmu)\| = P(\bfx - \bfmu)/\|\bfx - \bfmu\| = P \bfS(\bfx - \bfmu)$ for any orthogonal matrix $P$, and as a result the SCM is rotation equivariant too, in the sense that $\Sigma_S(P\bfX) = P \Sigma_S(\bfX) P^T$. This is not necessarily true in general if $P$ is replaced by any non-singular matrix. An affine equivariant version of the sample SCM is obtained as the solution $\hat \Sigma_T$ of the following equation:
$$ \hat \Sigma_T(\bfX) = \frac{p}{n} \sum_{i=1}^n \frac{(\bfx - \bfmu)(\bfx - \bfmu)^T}{(\bfx - \bfmu)^T \hat\Sigma_T(\bfX)^{-1} (\bfx - \bfmu)} $$
which turns out to be Tyler's M-estimator of scatter \citep{tyler87}. In this context, one should note that for scatter matrices, affine equivariance will mean any affine transformation on the original random variable $\bfX \mapsto \bfX^* = A\bfX + \bfb$ ($A$ non-singular, $\bfb \in \mathbb{R}^p$) being carried over to the covariance matrix estimate upto a scalar multiple: $\hat\Sigma_T(\bfX^*) = k. A \hat\Sigma_T(\bfX) A^T$ for some $k>0$.

\subsection{Data depth and outlyingness}
For any multivariate distribution $F = F_\bfX$ belonging to a set of distributions $\mathcal F$, the depth of a point $\bfx \in \mathbb{R}^p$, say $D(\bfx, F_\bfX)$ is any real-valued function that provides a 'center outward ordering' of $\bfx$ with respect to $F$ \citep{zuo00}. \cite{liu90} outlines the desirable properties of a statistical depth function:

\vspace{1em}
\noindent\textbf{(D1)} \textit{Affine invariance}: $D(A\bfx + \bfb, F_{A\bfX+\bfb}) = D(\bfx, F_\bfX)$;

\noindent\textbf{(D2)} \textit{Maximality at center}: $D(\bftheta, F_\bfX) = \sup_{\bfx\in \mathbb{R}^p} D(\bfx, F_\bfX)$ for $F_\bfX$ having center of symmetry $\bftheta$. This point is called the \textit{deepest point} of the distribution.;

\noindent\textbf{(D3)} \textit{Monotonicity with respect to deepest point}: $D(\bfx; F_\bfX) \leq D(\bftheta + a(\bfx - \bftheta), F_\bfX)$, $\bftheta$ being deepest point of $F_\bfX$.;

\noindent\textbf{(D4)} \textit{Vanishing at infinity}: $D(\bfx; F_\bfX) \rightarrow 0$ as $\|\bfx\| \rightarrow \infty $.
\vspace{1em}

In (D2) the types of symmetry considered can be central symmetry, angular symmetry and halfspace symmetry. Also for multimodal probability distributions, i.e. distributions with multiple local maxima in their probability density functions, properties (D2) and (D3) are actually restrictive towards the formulation of a reasonable depth function that captures the shape of the data cloud. In our derivations that follow, we replace these two by a weaker condition:

\vspace{1em}
\noindent\textbf{(D2*)} \textit{Existence of a maximal point}: The maximum depth over all distributions $F$ and points $\bfx$ is bounded above, i.e. $ \sup_{F_\bfX \in \mathcal F} \sup_{\bfx\in \mathbb{R}^p} D(\bfx, F_\bfX) < \infty $. We denote this point by $M_D(F_\bfX)$.
\vspace{1em}

A real-valued function measuring the outlyingness of a point with respect to the data cloud can be seen as the opposite of what data depth does. Indeed, such functions have been used to define several depth functions, for example simplicial depth, projection depth and $L_p$-depth. Keeping with the spirit of the utility of these functions we name them `htped': literally the reverse of `depth', and give a general definition of such functions as a transformation on any depth function:

\begin{Definition}
Given a random variable $\bfX$ following a probability distribution $F$, and a depth function $D(.,.)$, we define Htped of a point $\bfx$ as: $\tilde D(\bfx, F) = h(d_\bfx)$ as any function of the data depth $D(\bfx, F) = d_\bfx$ so that $h(d_\bfx)$ is bounded, monotonically decreasing in $d_\bfx$ and $\sup_\bfx \tilde D(\bfx, F) < \infty$.
\end{Definition}

For a fixed depth function, there are several choices of a corresponding htped. We develop our theory assuming a general htped function, but for the plots and simulations, fix our htped as $\tilde D(\bfx, F) = M_D(F) - D(\bfx, F)$, i.e. simply subtract the depth of a point from the maximum possible depth over all points in sample space.

We will be using the following 3 measures of data-depth to obtain our DCMs and compare their performances:

\begin{itemize}
\item \textbf{Halfspace depth} (HD) \citep{tukey75} is defined as the minimum probability of all halfspaces containing a point. In our notations,

$$ HD(\bfx, F)  = \inf_{\bfu \in \mathbb{R}^p; \bfu \neq \bf0} P(\bfu^T \bfX \geq \bfu^T \bfx) $$

\item \textbf{Mahalanobis depth} (MhD) \citep{LiuPareliusSingh99} is based on the Mahalanobis distance of $\bfx$ to $\bfmu$ with respect to $\Sigma$: $d_\Sigma(\bfx, \bfmu) = \sqrt{(\bfx - \bfmu)^T \Sigma^{-1} (\bfx - \bfmu)}$. It is defined as
$$ MhD(\bfX, F) = \frac{1}{1 + d^2_\Sigma (\bfx - \bfmu)} $$
note here that $d_\Sigma(\bfx,\bfmu)$ can be seen as a valid htped function of $\bfx$ with respect to $F$.

\item \textbf{Projection depth} (PD) \citep{zuo03} is another depth function based on an outlyingness function. Here that function is
$$ O(\bfx, F) = \sup_{\| \bfu \| = 1} \frac{| \bfu^T\bfx - m(\bfu^T\bfX)|}{s(\bfu^T\bfX)} $$
where $m$ and $s$ are some univariate measures location and scale, respectively. Given this the depth at $\bfx$ is defined as $PD(\bfx, F) = 1/(1+O(\bfx, F))$.
\end{itemize}

Computation-wise, MhD is easy to calculate since the sample mean and covariance matrix are generally used as estimates of $\mu$ and $\Sigma$, respectively. However this makes MhD less robust with respect to outliers. PD is generally approximated by taking maximum over a number of random projections. There have been several approaches for calculating HD. A recent unpublished paper \citep{rainerArxiv} provides a general algorithm that computes exact HD in $O(n^{p-1}\log n)$ time. In this paper, we shall use inbuilt functions in the R package \texttt{fda.usc} for calculating the above depth functions.

\section{Depth-based rank covariance matrix} \label{section:dcmSection}

Consider a vector-valued random variable $\bfX \in \mathbb{R}^p$ . Data depth is as much a property of the random variable as it is of the underlying distribution, so for ease of notation while working with transformed random variables, from now on we shall be using $D_\bfX(\bfx) = D(\bfx, F)$ to denote the depth of a point $\bfx$. Now, given a depth function $D_{\bfX}(\bfx)$ (equivalently, an htped function $\tilde D_\bfX(\bfx) = \tilde D(\bfx, F)$), transform the original random variable as: $\tilde \bfx = \tilde D_\bfX(\bfx) \bfS(\bfx - \bfmu)$, $\bfS(.)$ being the spatial sign functional. The transformed random variable $\tilde \bfX$ can be seen as the multivariate rank corresponding to $\bfX$ (e.g. \cite{serfling2006}). The notion of multivariate ranks goes back to \cite{PuriSenBook}, where they take the vector consisting of marginal univariate ranks as multivariate rank vector. Subsequent definitions of multivariate ranks were proposed by \cite{MottonenOja95,HallinPaindaveine02} and \cite{Chernozhukov14}. Compared to these formulations, our definition of multivariate ranks works for any general depth function, and provides an intuitive extension to any spatial sign-based methodology.

\begin{figure}[t]
	\captionsetup{singlelinecheck=off}
	\centering
		\includegraphics[height=6cm]{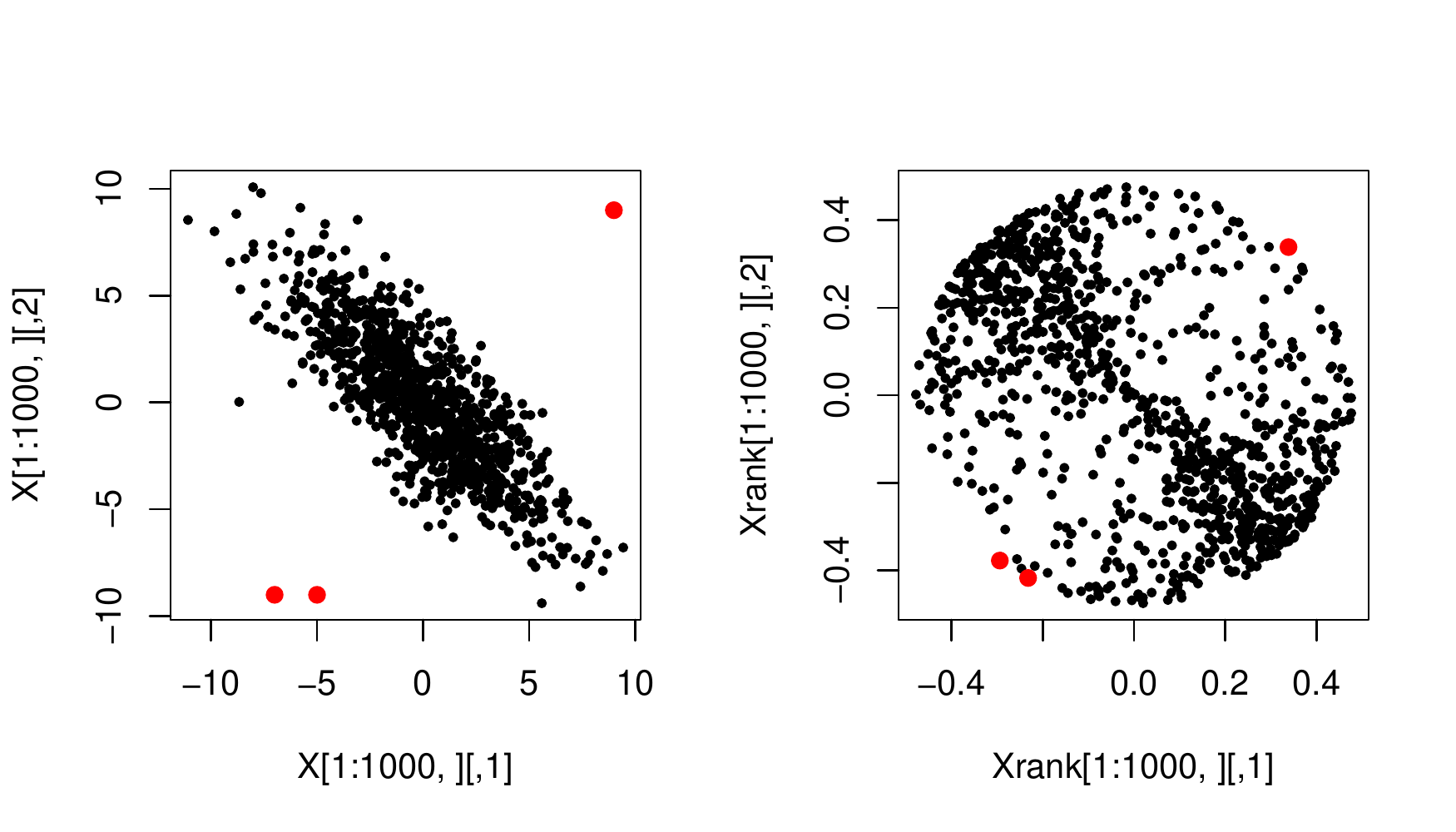}
	\caption{(Left) 1000 points randomly drawn from $\mathcal N_2\left((0,0)^T, \left(\protect\begin{smallmatrix} 5 & -4 \\ -4 & 5 \protect\end{smallmatrix}\right)\right) $ and (Right) their multivariate ranks based on halfspace depth}
	\label{fig:rankplot}
\end{figure}

Figure \ref{fig:rankplot} gives an idea of how the multivariate rank vector $\tilde \bfX$ is distributed when $\bfX$ has a bivariate normal distribution. Compared to the spatial sign, which are distributed on the surface of $p$-dimensional unit ball centered at $\bfmu$, these spatial ranks have the same direction as original data and reside \textit{inside} the $p$-dimensional ball around $\bfmu$ that has radius $M_D(F)$ (which, for the case of halfspace depth, equals 0.5). Any outlying samples situated far away from the data cloud (represented by red points in the figure) are mapped close to the boundary of the $p$-dimensional ball after the rank transformation.

Now consider the spectral decomposition for the covariance matrix of $F$: $\Sigma = \Gamma\Lambda\Gamma^T$, $\Gamma$ being orthogonal and $\Lambda$ diagonal with positive diagonal elements. Also normalize the original random variable as $\bfz = \Gamma^T\Lambda^{-1/2} (\bfx - \bfmu)$. In this setup, we can represent the transformed random variable as
\begin{eqnarray}
\tilde \bfx &=& \tilde D_{\bfX} (\bfx) \bfS(\bfX - \bfmu) \notag \\
&=& \tilde D_{\Gamma\Lambda^{1/2}\bfZ + \bfmu} (\Gamma\Lambda^{1/2} \bfz + \bfmu). \bfS(\Gamma\Lambda^{1/2} \bfz) \notag \\
&=& \Gamma \tilde D_{\bfZ}(\bfz) \bfS(\Lambda^{1/2}\bfz) \notag \\
&=& \Gamma \Lambda^{1/2} \tilde D_{\bfZ}(\bfz) \bfS(\bfz) \frac{\| \bfz \|}{\|\Lambda^{1/2} \bfz\|}
\label{equation:rankdecomp}
\end{eqnarray}
%

$\tilde D_\bfZ(\bfz)$ is an even function in $\bfz$ because of affine invariance, as is $\| \bfz \| / \|\Lambda^{1/2} \bfz \|$. Since $\bfS(\bfz)$ is odd in $\bfz$ for circularly symmetric $\bfz$, it follows that $E(\tilde \bfX) = \bf0$, and consequently we obtain an expression for the covariance matrix of $\tilde \bfX$:

\begin{Theorem} \label{Theorem:covform}
Let the random variable $\bfX \in \mathbb{R}^p$ follow an elliptical distribution with center $\bfmu$ and covariance matrix $\Sigma = \Gamma\Lambda\Gamma^T$, its spectral decomposition. Then, given a depth function $D_\bfX(.)$ the covariance matrix of the transformed random variable $\tilde\bfX$ is
\begin{equation} \label{equation:covformEq1}
Cov(\tilde \bfX) = \Gamma \Lambda_{D,S} \Gamma^T,\quad\mbox{with}\quad \Lambda_{D,S} = \mathbb E_\bfz \left[ (\tilde D_\bfZ(\bfz))^2 \frac{\Lambda^{1/2} \bfz \bfz^T \Lambda^{1/2}}{\bfz^T \Lambda \bfz} \right]
\end{equation}
where $\bfz = (z_1,...,z_p)^T \sim N({\bf 0}, I_p)$, so that $\Lambda_{D,S}$ a diagonal matrix with diagonal entries
$$ \lambda_{D,S,i} = \mathbb E_\bfZ \left[ \frac{(\tilde D_\bfZ(\bfz))^2 \lambda_i z_i^2}{\sum_{j=1}^p \lambda_j z_j^2} \right] $$
\end{Theorem}

The matrix of eigenvectors of the covariance matrix, $\Gamma$, remains unchanged in the transformation $\bfX \rightarrow \tilde \bfX$. As a result, the multivariate rank vectors can be used for robust principal component analysis, which will be outlined in the following sections. However, as one can see in the above expression, the diagonal entries of $\Lambda_{D,S}$ are not same as those of $\Lambda$, i.e. the actual eigenvalues. This is the reason for lack of affine equavariance of the DCM. Following the case of multivariate sign covariance matrices \citep{taskinen12} one can get back the shape components, i.e. original \textit{standardized} eigenvalues $\Lambda^*$ from $\Lambda_{D,S}$ by an iterative algorithm:

\begin{enumerate}
\item Set $k=0$, and start with an initial value $\Lambda^{*(0)}$.

\item Calculate the next iterate
$$ \Lambda^{*(k+1)} = \left[ \mathbb E_\bfz \left( \frac{(\tilde D_\bfZ(\bfz))^2 \bfz \bfz^T}{\bfz^T \Lambda^{*(k)} \bfz} \right) \right]^{-1} \Lambda_{D,S} $$
and standardize its eigenvalues:
$$ \Lambda^{*(k+1)} = \frac{\Lambda^{*(k+1)}}{\text{det} (\Lambda^{*(k+1)})^{1/p}} $$
\item Stop if convergence criterion is satisfied. Otherwise set $k \rightarrow k+1$ ad go to step 2.
\end{enumerate}

Unlike sign covariance matrices and symmetrized sign covariance matrices \citep{dumbgen98}, however, attempting to derive an affine equivariant counterpart (as opposed to only orthogonal equivariance) of the DCM through an iterative approach analogous to \cite{tyler87} will not result in anything new. This is because Tyler's scatter matrix $\Sigma_T$ is defined as the implicit solution to the following equation:
\begin{equation} \label{equation:tylerEq}
\Sigma_T = \mathbb E\left[ \frac{(\bfx - \bfmu) (\bfx - \bfmu)^T}{(\bfx - \bfmu)^T \Sigma_T^{-1} (\bfx - \bfmu)} \right]
\end{equation}
and simply replacing $\bfx$ by its multivariate rank counterpart $\tilde\bfx$ will not change the estimate $\Sigma_T$ as $\bfx$ and $\tilde \bfx$ have the same directions. Instead we consider a depth-weighted version of Tyler's scatter matrix (i.e. weights $(\tilde D_\bfX(\bfx))^2$ in right side of (\ref{equation:tylerEq})) in the simulations in Section \ref{section:simSection}. The simulations show that it has slightly better finite-sample efficiency than $\Sigma_T$ but has same asymptotic performance. We conjecture that its concentration properties can be obtained by taking an approach similar to \cite{soloveychik14}.

\section{Asymptotic results}\label{section:sec4}
\subsection{The sample DCM}
Let us now consider $n$ iid random draws from our elliptic distribution $F$, say $\bfX_1,...,\bfX_n$. For ease of notation, denote $SS(\bfx; \bfmu) = \bfS(\bfx - \bfmu) \bfS(\bfx - \bfmu)^T$. Then, given the depth function and known location center $\bfmu$, one can show that the vectorized form of $\sqrt n$-times the sample DCM: $\sum_{i=1}^n (\tilde D_\bfX(\bfx_i))^2SS(\bfx_i; \bfmu) /\sqrt{n}$  has an asymptotic multivariate normal distribution with mean $\sqrt n.vec(E[( (\tilde D_\bfX (\bfX))^2 SS(\bfx; \bfmu)])$ and a certain covariance matrix by straightforward application of the central limit theorem (CLT). But in practice the population depth function $D_\bfX(\bfx) = D(\bfx, F)$ is estimated by the depth function based on the empirical distribution function, $F_n$. Denote this sample depth by $D^n_\bfX (\bfx) = D(\bfx, F_n)$. Here we make the following assumption regarding its relation to $D_\bfX(\bfx)$:

\vspace{1em}
\noindent\textbf{(D5)} \textit{Uniform convergence}: $\sup_{\bfx \in \mathbb R^p} | D^n_\bfX (\bfx) - D_\bfX (\bfx) | \rightarrow 0$ as $n \rightarrow \infty $.
\vspace{1em}

The assumption that empirical depths converge uniformly at all points $\bfx$ to their population versions holds under very mild conditions for several well known depth functions: for example projection depth \citep{zuo03} and simplicial depth \citep{Dumbgen92}. One also needs to replace the known location parameter $\bfmu$ by some estimator $\hat\bfmu_n$. Examples of robust estimators of location that are relevant here include the spatial median \citep{haldane48,brown83}, Oja median \citep{oja83}, projection median \citep{zuo03} etc. Now, given $D^n_\bfX(.)$ and $\hat \bfmu_n$, to plug them into the sample DCM and still go through with the CLT we need the following result:

\begin{Lemma} \label{Lemma:lemma1}
Consider a random variable $\bfX \in \mathbb{R}^p$ having a continuous and symmetric distribution with location center $\bfmu$ such that $E\|\bfx - \bfmu \|^{-3/2} < \infty$. Given $n$ random samples from this distribution, suppose $\hat\bfmu_n$ is an estimator of $\bfmu$ so that $\sqrt n (\hat\bfmu_n - \bfmu) = O_P(1) $. Then with the above notations, and given the assumption (D5) we have
$$ \sqrt n \left[
\frac{1}{n} \sum_{i=1}^n (\tilde D^n_\bfX (\bfx_i))^2 SS(\bfx_i; \hat\bfmu_n) -
\frac{1}{n} \sum_{i=1}^n (\tilde D_\bfX (\bfx_i))^2 SS(\bfx_i; \bfmu) \right]
\stackrel{P}{\rightarrow} 0 $$
\end{Lemma}

We are now in a position to state the result for consistency of the sample DCM:

\begin{Theorem} \label{Theorem:rootn}
Consider $n$ iid samples from the distribution in Lemma \ref{Lemma:lemma1}. Then, given a depth function $D_\bfX(.)$ and an estimate of center $\hat\bfmu_n$ so that $\sqrt n(\hat \bfmu_n - \bfmu) = O_P(1)$,
$$ \sqrt n \left[ vec\left\{ \frac{1}{n} \sum_{i=1}^n (\tilde D^n_\bfX (\bfx_i))^2 SS(\bfx_i; \hat\bfmu_n) \right\} - E \left[ vec\left\{ (\tilde D_\bfX (\bfx))^2 SS(\bfx; \bfmu) \right\} \right] \right]
\stackrel{D}{\rightarrow}
N_{p^2} ({\bf 0}, V_{D,S}(F)) $$
$$ \text{with } V_{D,S}(F) = Var \left[vec \left\{ (\tilde D_\bfX (\bfx))^2 SS(\bfx; \bfmu) \right\} \right] $$
\end{Theorem}

In case $F$ is elliptical, an elaborate form of the covariance matrix $V_{D,S}(F)$ explicitly specifying each of its elements (more directly those of its $\Gamma^T$-rotated version) can be obtained, which is given in Appendix \ref{section:appA}. This form is useful when deriving limiting distributions of eigenvectors and eigenvalues of the sample DCM.

\subsection{Eigenvectors and eigenvalues} Since we are mainly interested in using the DCM for a robust version of principal components analysis, from now on we assume that the eigenvalues of $\Sigma$ are distinct: $\lambda_1 > \lambda_2 > ... > \lambda_p$ to obtain asymptotic distributions of principal components. In the case of eigenvalues with larger than 1 multiplicities, the limiting distributions of eigenprojection matrices can be obtained analogous to those of the sign covariance matrix \citep{magyar14}.

\paragraph{}We now derive the asymptotic joint distributions of eigenvectors and eigenvalues of the sample DCM. The following result allows us to get these, provided we know the limiting distribution of the sample DCM itself:

\begin{Theorem} \label{Theorem:decomp} \citep{taskinen12}
Let $F_\Lambda$ be an elliptical distribution with a diagonal covariance matrix $\Lambda$, and $\hat C$ be any positive definite symmetric $p \times p$ matrix such that at $F_\Lambda$ the limiting distribution of $\sqrt n vec(\hat C - \Lambda)$ is a $p^2$-variate (singular) normal distribution with mean zero. Write the spectral decomposition of $\hat C$ as $\hat C = \hat P \hat\Lambda \hat P^T$. Then the limiting distributions of $\sqrt n vec(\hat P - I_p)$ and $\sqrt n vec(\hat\Lambda - \Lambda)$ are multivariate (singular) normal and
\begin{equation} \label{equation:decompEq}
\sqrt n vec (\hat C - \Lambda)  = \left[ (\Lambda \otimes I_p) - (I_p \otimes \Lambda) \right] \sqrt n vec (\hat P - I_p) + \sqrt n vec (\hat\Lambda - \Lambda) + o_P(1)
\end{equation}
\end{Theorem}

The first matrix picks only off-diagonal elements of the LHS and the second one only diagonal elements. We shall now use this as well as the form of the asymptotic covariance matrix of the vec of sample DCM, i.e. $V_{D,S}(F)$ to obtain limiting variance and covariances of eigenvalues and eigenvectors.

\begin{Corollary} \label{Corollary:eigendist}
Consider the sample DCM $ \hat S^D(F) = \sum_{i=1}^n (\tilde D^n_\bfX (\bfx_i))^2 SS(\bfx_i; {\bf \hat\bfmu_n})/n $ and its spectral decomposition $\hat S^D(F) = \hat\Gamma_D \hat\Lambda_D \hat\Gamma_D^T $. Then the matrices $G = \sqrt n (\hat\Gamma_D - \Gamma) $ and $L = \sqrt n (\hat\Lambda_D - \Lambda_{D,S}) $ have independent distributions. The random variable $vec(G)$ asymptotically has a $p^2$-variate normal distribution with mean ${\bf 0}_{p^2}$, and the asymptotic variance and covariance of different columns of $G = (\bfg_1,...,\bfg_p)$ are as follows:
\begin{equation} \label{equation:DevEq}
AVar(\bfg_i) = \sum_{k=1; k \neq i}^p \frac{1}{(\lambda_{D,s,k} - \lambda_{D,S,i})^2} E \left[ \frac{(\tilde D_\bfZ (\bfz))^4 \lambda_i \lambda_k z_i^2 z_k^2}{(\bfz^T \Lambda \bfz)^2} \right] \bfgamma_k \bfgamma_k^T
\end{equation}
\begin{equation}
ACov(\bfg_i, \bfg_j) = - \frac{1}{(\lambda_{D,s,i} - \lambda_{D,S,j})^2} E \left[ \frac{(\tilde D_\bfZ (\bfz))^4 \lambda_i \lambda_j z_i^2 z_j^2}{(\bfz^T \Lambda \bfz)^2} \right] \bfgamma_j \bfgamma_i^T; \quad i \neq j
\end{equation}
where $\Gamma = (\bfgamma_1,...,\bfgamma_p)$. The vector consisting of diagonal elements of $L$, say $\bfl = (\l_1,...,\l_p)^T$ asymptotically has a $p$-variate normal distribution with mean ${\bf 0}_p$ and variance-covariance elements:
\begin{eqnarray}
AVar(l_i) &=& E \left[ \frac{(\tilde D_\bfZ (\bfz))^4 \lambda_i^2 z_i^4}{(\bfz^T \Lambda \bfz)^2} \right] - \lambda_{D,S,i}^2\\
ACov(l_i, l_j) &=& E \left[ \frac{(\tilde D_\bfZ (\bfz))^4 \lambda_i \lambda_j z_i^2 z_j^2}{(\bfz^T \Lambda \bfz)^2} \right] - \lambda_{D,S,i} \lambda_{D,S,j}; \quad i \neq j
\end{eqnarray}
\end{Corollary}

\section{Robustness and efficiency properties} \label{section:simSection}
In this section, we first obtain the influence functions of the DCM as well as its eigenvectors and eigenvalues, which are essential to understand how much influence a sample point, especially an infinitesimal contamination, has on any functional on the distribution \citep{hampel}. We also derive the asymptotic efficiencies of individual principal components with respect to those of the original covariance matrix and sign covariance matrix. Unlike affine equivariant estimators of shape, the Asymptotic Relative Efficiency (ARE) of eigenvectors (with respect to any other affine equivariant estimator) can not be simplified as a ratio of two scalar quantities dependent on only the distribution of $\| \bfz \|$ (e.g. \cite{taskinen12,ollilia03}). Finite sample efficiency of the DCM estimates with respect to infinitesimal contamination and heavy-tailed distributions shall also be demonstrated by a simulation study.

\subsection{Influence functions}
Given any probability distribution $F$, the influence function of any point $\bfx_0$ in the sample space $\mathcal{X}$ for some functional $T(F)$ on the distribution is defined as
$$ IF(\bfx_0; T,F) = \lim_{\epsilon \rightarrow 0} \frac{1}{\epsilon} (T(F_\epsilon) - T(F)) $$
where $F_\epsilon$ is $F$ with an additional mass of $\epsilon$ at $\bfx_0$, i.e. $F_\epsilon = (1-\epsilon)F + \epsilon \Delta_{\bfx_0}$; $\Delta_{\bfx_0}$ being the distribution with point mass at $\bfx_0$. When $T(F) = E_F g$ for some $F$-integrable function $g$, $IF(\bfx_0; T,F) = g(\bfx_0) - T(F)$. It now follows that for the DCM,
$$ IF(\bfx_0; Cov(\tilde \bfX), F) = (\tilde D_\bfX(\bfx_0))^2 SS(\bfx_0; \bfmu) - Cov(\tilde \bfX) $$

\begin{figure}[]
	\centering
		\includegraphics[width=12cm]{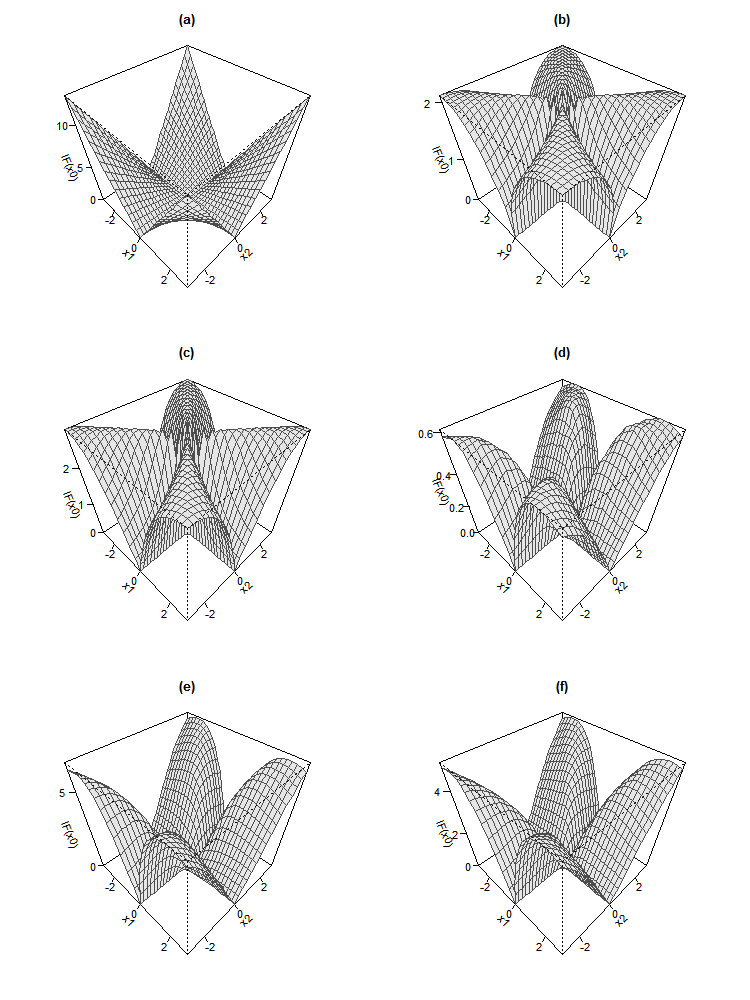}
	\caption{Plot of the norm of influence function for first eigenvector of (a) sample covariance matrix, (b) SCM, (c) Tyler's scatter matrix and DCMs for (d) Halfspace depth, (e) Mahalanobis depth, (f) Projection depth for a bivariate normal distribution with $\bfmu = {\bf 0}, \Sigma = \diag(2,1)$}
	\label{fig:IFnorm}
\end{figure}

Following \cite{croux00}, we now get the influence function of the $i^\text{th}$ column of $\hat\Gamma_D = (\hat\bfgamma_{D,1},...,\hat\bfgamma_{D,p}); i = 1,...,p$:
\begin{eqnarray}
IF(\bfx_0; \hat\bfgamma_{D,i}, F) &=& \sum_{k=1; k \neq i}^p \frac{1}{\lambda_{D,S,i} - \lambda_{D,S,k}} \left\{ \bfgamma^T_k IF(\bfx_0; Cov(\tilde \bfX), \bfgamma_i) \right\} \bfgamma_k \notag \\
&=& \sum_{k=1; k \neq i}^p \frac{1}{\lambda_{D,S,i} - \lambda_{D,S,k}} \left\{ \bfgamma^T_k (\tilde D_\bfX(\bfx_0))^2 SS(\bfx_0; \bfmu)\bfgamma_i - \lambda_{D,S,i}\bfgamma_k^T\bfgamma_i \right\} \bfgamma_k \notag \\
&=& \sum_{k=1; k \neq i}^p \frac{\sqrt{\lambda_i \lambda_k} z_{0i} z_{0k}}{\lambda_{D,S,i} - \lambda_{D,S,k}}. \frac{(\tilde D_\bfZ(\bfz_0))^2 }{\bfz_0^T \Lambda \bfz_0} \bfgamma_k
\end{eqnarray}
where $\Gamma^T \Lambda^{-1/2} (\bfx_0 - \bfmu) = \bfz_0 = (z_{01},...,z_{0p})^T$. Clearly this influence function will be bounded, which indicates good robustness properties of principal components. Moreover, since the htped function takes small values for points close to the center of the distribution, it does not suffer from the inlier effect that is typical of the SCM and Tyler's shape matrix. The influence function for the $i^\text{th}$ eigenvector estimates of these two matrices (say $\hat\bfgamma_{S,i}$ and $\hat\bfgamma_{T,i}$, respectively) are as follows:
\begin{eqnarray*}
\quad IF(\bfx_0; \hat \bfgamma_{S,i}, F) &=& \sum_{k=1; k \neq i}^p \frac{\sqrt{\lambda_i \lambda_k}}{\lambda_{S,i} - \lambda_{S,k}}. \frac{z_{0i} z_{0k}}{\bfz_0^T \Lambda \bfz_0} \bfgamma_k, \text{ with } \lambda_{S,i} = E_\bfZ \left( \frac{\lambda_i z_i^2}{\sum_{j=1}^p \lambda_j z_j^2} \right) \quad \\
IF(\bfx_0; \hat \bfgamma_{T,i}, F) &=& (p+2) \sum_{k=1; k \neq i}^p \frac{\sqrt{\lambda_i \lambda_k}}{\lambda_i - \lambda_k}. \frac{z_{0i} z_{0k}}{\bfz_0^T \bfz_0} \bfgamma_k \quad 
\end{eqnarray*}
for $i = 1,2$. In Figure \ref{fig:IFnorm} we consider first eigenvectors of different scatter estimates for the $\mathcal{N}_2((0,0)^T, \diag(2,1))$ and plot norms of these influence functions for different values of $\bfx_0$. The plots for SCM and Tyler's shape matrix demonstrate the 'inlier effect', i.e. points close to symmetry center and the center itself having high influence. The influence function for the sample covariance matrix is obtained by replacing $(p+2)$ by $\| \bfz_0 \|^2$ in the expression of $IF(\bfx_0; \hat \bfgamma_{T,i}, F)$ above, hence is unbounded and the corresponding eigenvector estimators are not robust. In comparison, all three DCMs considered here have a bounded influence function as well as small values of the influence function at 'deep' points.

\subsection{Asymptotic and finite-sample efficiencies}
Suppose $\hat\Sigma$ is a $\sqrt n$-consistent estimator of the population covariance matrix $\Sigma$, which permits a spectral decomposition $ \hat\Sigma = \hat\Gamma \hat\Lambda \hat\Gamma^T $, where $\hat\Gamma = (\hat\bfgamma_1,...,\hat\bfgamma_p)$. Then the asymptotic variance of the eigenvectors are (see Theorem 13.5.1 in \cite{anderson})
\begin{equation} \label{equation:covevEq}
AVar(\sqrt n\hat \bfgamma_i) = \sum_{k=1; k \neq i}^p \frac{\lambda_i \lambda_k}{(\lambda_i - \lambda_k)^2} \bfgamma_k \bfgamma_k^T
\end{equation}
The asymptotic relative efficiencies of eigenvectors from the sample DCM with respect to the sample covariance matrix can now be derived using (\ref{equation:covevEq}) above and (\ref{equation:DevEq}) from Corollary \ref{Corollary:eigendist}:
\begin{eqnarray*}
ARE(\hat\bfgamma^D_i, \hat\bfgamma_i; F) &=& \frac{\Tr( AVar(\sqrt n\hat \bfgamma_i))}{\Tr( AVar(\sqrt n\hat \bfgamma^D_i))}\\
&=& \left[\sum_{k=1; k \neq i}^p \frac{\lambda_i \lambda_k}{(\lambda_i - \lambda_k)^2} \right] \left[ \sum_{k=1; k \neq i}^p \frac{\lambda_i \lambda_k }{(\lambda_{D,s,i} - \lambda_{D,S,k})^2} E \left( \frac{(\tilde D_\bfZ (\bfz))^4 z_i^2 z_k^2}{(\bfz^T \Lambda \bfz)^2} \right) \right]^{-1}
\end{eqnarray*}

For 2 dimensions, this expression can be somewhat simplified. Suppose the two eigenvalues are $\lambda$ and $\rho\lambda$. In that case the eigenvalues of the DCM are
$$ \lambda_{D,S,1} = E \left( \frac{(\tilde D_\bfZ(\bfz))^2 z_1^2}{z_1^2 + \rho z_2^2} \right), \quad
\lambda_{D,S,2} = E \left( \frac{(\tilde D_\bfZ(\bfz))^2 \rho z_2^2}{z_1^2 + \rho z_2^2} \right) $$
and by simple algebra we get
$$ ARE(\hat\bfgamma^D_1, \hat\bfgamma_1; F) = ARE(\hat\bfgamma^D_2, \hat\bfgamma_2; F) = \frac{1}{(1-\rho)^2} \frac{\left[ E \left( \frac{(\tilde D_\bfZ(\bfz))^2 (z_1^2 - \rho z_2^2)}{(z_1^2 + \rho z_2^2)} \right) \right]^2 }{E \left( \frac{(\tilde D_\bfZ(\bfz))^4 z_1^2z_2^2}{(z_1^2 + \rho z_2^2)^2}\right)} $$

For $\rho=0.5$, Table \ref{table:AREtable} below considers 6 different elliptic distributions (namely, bivariate $t$ with df = $5,6,10,15,25$ and bivariate normal), and summarizes the ARE for first eigenvector of the SCM, Tyler's scatter matrix and DCM for 3 choices of depth function (HSD-CM, MhD-CM, PD-CM: columns 3-5), as well as their depth-weighted versions mentioned at the end of section \ref{section:dcmSection} (HSD-wCM, MhD-wCM, PD-wCM: columns 6-8). The SCM and Tyler's M-estimator perform better than the sample covariance matrix only for a bivariate $t$-distribution with df = 5. Estimates based on depth-based covariance matrices have much better performances for all distributions, and continue to be competitive of the sample covariance matrix estimates as the base distribution approaches normality, especially those based on projection depth. Interestingly, depth-weighted versions seem to have better performances than their corresponding DCMs, more so for heavy-tailed distributions.

\begin{table}
\begin{footnotesize}
    \begin{tabular}{c|cc|ccc|ccc}
     \hline
                       & SCM  & Tyler & HSD-CM & MhD-CM & PD-CM & HSD-wCM & MhD-wCM & PD-wCM \\ \hline
    Bivariate $t_5$    & 1.46 & 1.50  & 1.97   & 1.54   & 1.99  & 2.09    & 1.57    & 2.11   \\
    Bivariate $t_6$    & 0.97 & 1.00  & 1.37   & 1.12   & 1.40  & 1.45    & 1.19    & 1.49   \\
    Bivariate $t_{10}$ & 0.65 & 0.67  & 0.96   & 0.88   & 1.02  & 1.02    & 0.94    & 1.10   \\
    Bivariate $t_{15}$ & 0.57 & 0.59  & 0.87   & 0.82   & 0.94  & 0.92    & 0.87    & 1.00   \\
    Bivariate $t_{25}$ & 0.54 & 0.55  & 0.79   & 0.75   & 0.88  & 0.85    & 0.81    & 0.95   \\
    BVN   & 0.49 & 0.50  & 0.73   & 0.71   & 0.82  & 0.77    & 0.75    & 0.88   \\ \hline
    \end{tabular}
\end{footnotesize}
\caption{Asymptotic efficiencies relative to sample covariance matrix for $p=2$}
\label{table:AREtable}
\end{table}

We now obtain finite sample efficiencies of the three DCMs as well as their depth-weighted affine equivariant counterparts by a simulation study, and compare them with the same from the SCM and Tyler's scatter matrix. We consider the same 6 elliptical distributions considered in ARE calculations above, and from every distribution draw 10000 samples each for sample sizes $n = 20, 50, 100, 300, 500$. All distributions are centered at ${\bf 0}_p$, and have covariance matrix $\Sigma = \diag(p,p-1,...1)$. We consider 3 choices of $p$: 2, 3 and 4.

We use the concept of principal angles \citep{miao92} to find out error estimates for the first eigenvector of a scatter matrix. In our case, the first eigenvector will be
$$ \bfgamma_1 = (1,\overbrace{0,...,0}^{p-1})^T $$
For an estimate of the eigenvector, say $\hat\bfgamma_1$, error in prediction is measured by the smallest angle between the two lines, i.e. $ \cos^{-1} | \bfgamma_1^T \hat\bfgamma_1 | $. A smaller absolute value of this angle is equivalent to better prediction. We repeat this 10000 times and calculate the \textbf{Mean Squared Prediction Angle}:
$$ MSPA(\hat \bfgamma_1) = \frac{1}{10000} \sum_{m=1}^{10000} \left( \cos^{-1} \left|\bfgamma_1^T \hat\bfgamma^{(m)}_1 \right| \right)^2 $$
Finally, the finite sample efficiency of some eigenvector estimate $\hat\bfgamma^E_1$ relative to that obtained from the sample covariance matrix, say $\hat\bfgamma^{Cov}_1$ is obtained as:
$$ FSE(\hat\bfgamma^E_1, \hat\bfgamma^{Cov}_1) = \frac{MSPA(\hat\bfgamma^{Cov}_1)}{MSPA(\hat\bfgamma^E_1)} $$

Tables \ref{table:FSEtable2}, \ref{table:FSEtable3} and \ref{table:FSEtable4} give FSE values for $p=2,3,4$, respectively. In general, all the efficiencies increase as the dimension $p$ goes up. DCM-based estimators (columns 3-5 in each table) outperform SCM and Tyler's scatter matrix, and among the 3 depths considered, projection depth seems to give the best results. Its finite sample performances are better than Tyler's and Huber's M-estimators of scatter as well as their symmetrized counterparts (see Table 4 in \cite{sirkia07}, and quite close to the affine equivariant spatial sign covariance matrix (see Table 2 in \cite{ollilia03}) For $p=2$, $n=300, 500$ the first 5 columns of Table \ref{table:FSEtable2} approximate the asymptotic efficiencies in Table \ref{table:AREtable} well, except for the multivariate $t$-distribution with df = 5. Finally, the depth-weighted iterated versions of these 3 SCMs (columns 6-8 in each table) seem to further better the performance of their corresponding orthogonal equivariant counterparts.

\begin{table}
\begin{footnotesize}
    \begin{tabular}{c|cc|ccc|ccc}
    \hline
    $F$ = Bivariate $t_5$    & SCM  & Tyler & HSD-CM & MhD-CM & PD-CM & HSD-wCM & MhD-wCM & PD-wCM \\ \hline
    $n$=20                   & 0.80 & 0.83  & 0.95   & 0.95   & 0.89  & 1.00    & 0.96    & 0.89   \\
    $n$=50                   & 0.86 & 0.90  & 1.25   & 1.10   & 1.21  & 1.32    & 1.13    & 1.25   \\
    $n$=100                  & 1.02 & 1.04  & 1.58   & 1.20   & 1.54  & 1.67    & 1.24    & 1.63   \\
    $n$=300                  & 1.24 & 1.28  & 1.81   & 1.36   & 1.82  & 1.93    & 1.44    & 1.95   \\
    $n$=500                  & 1.25 & 1.29  & 1.80   & 1.33   & 1.84  & 1.91    & 1.39    & 1.97   \\ \hline
    $F$ = Bivariate $t_6$    & SCM  & Tyler & HSD-CM & MhD-CM & PD-CM & HSD-wCM & MhD-wCM & PD-wCM \\ \hline
    $n$=20                   & 0.77 & 0.79  & 0.92   & 0.92   & 0.86  & 0.96    & 0.92    & 0.85   \\
    $n$=50                   & 0.76 & 0.78  & 1.11   & 1.00   & 1.08  & 1.17    & 1.03    & 1.13   \\
    $n$=100                  & 0.78 & 0.79  & 1.27   & 1.06   & 1.33  & 1.35    & 1.11    & 1.41   \\
    $n$=300                  & 0.88 & 0.91  & 1.29   & 1.09   & 1.35  & 1.38    & 1.15    & 1.45   \\
    $n$=500                  & 0.93 & 0.96  & 1.37   & 1.13   & 1.40  & 1.44    & 1.19    & 1.48   \\ \hline
    $F$ = Bivariate $t_{10}$ & SCM  & Tyler & HSD-CM & MhD-CM & PD-CM & HSD-wCM & MhD-wCM & PD-wCM \\ \hline
    $n$=20                   & 0.70 & 0.72  & 0.83   & 0.84   & 0.77  & 0.89    & 0.87    & 0.79   \\
    $n$=50                   & 0.58 & 0.60  & 0.90   & 0.84   & 0.86  & 0.95    & 0.88    & 0.91   \\
    $n$=100                  & 0.57 & 0.59  & 0.92   & 0.87   & 0.97  & 0.98    & 0.90    & 1.03   \\
    $n$=300                  & 0.62 & 0.64  & 0.93   & 0.85   & 0.99  & 0.99    & 0.91    & 1.06   \\
    $n$=500                  & 0.62 & 0.65  & 0.93   & 0.86   & 1.00  & 1.00    & 0.92    & 1.08   \\ \hline
    $F$ = Bivariate $t_{15}$ & SCM  & Tyler & HSD-CM & MhD-CM & PD-CM & HSD-wCM & MhD-wCM & PD-wCM \\ \hline
    $n$=20                   & 0.63 & 0.66  & 0.76   & 0.78   & 0.72  & 0.81    & 0.81    & 0.73   \\
    $n$=50                   & 0.52 & 0.52  & 0.79   & 0.75   & 0.80  & 0.84    & 0.79    & 0.85   \\
    $n$=100                  & 0.51 & 0.52  & 0.83   & 0.77   & 0.88  & 0.88    & 0.81    & 0.94   \\
    $n$=300                  & 0.55 & 0.56  & 0.84   & 0.79   & 0.91  & 0.89    & 0.84    & 0.98   \\
    $n$=500                  & 0.56 & 0.59  & 0.85   & 0.80   & 0.93  & 0.91    & 0.86    & 0.99   \\ \hline
    $F$ = Bivariate $t_{25}$ & SCM  & Tyler & HSD-CM & MhD-CM & PD-CM & HSD-wCM & MhD-wCM & PD-wCM \\ \hline
    $n$=20                   & 0.63 & 0.65  & 0.77   & 0.79   & 0.74  & 0.80    & 0.81    & 0.74   \\
    $n$=50                   & 0.49 & 0.50  & 0.73   & 0.71   & 0.76  & 0.78    & 0.75    & 0.80   \\
    $n$=100                  & 0.45 & 0.46  & 0.73   & 0.69   & 0.81  & 0.78    & 0.73    & 0.87   \\
    $n$=300                  & 0.51 & 0.52  & 0.78   & 0.75   & 0.87  & 0.83    & 0.79    & 0.94   \\
    $n$=500                  & 0.53 & 0.55  & 0.79   & 0.75   & 0.87  & 0.84    & 0.80    & 0.94   \\ \hline
    $F$ = BVN                & SCM  & Tyler & HSD-CM & MhD-CM & PD-CM & HSD-wCM & MhD-wCM & PD-wCM \\ \hline
    $n$=20                   & 0.56 & 0.60  & 0.69   & 0.71   & 0.67  & 0.73    & 0.74    & 0.68   \\
    $n$=50                   & 0.42 & 0.43  & 0.66   & 0.66   & 0.70  & 0.71    & 0.69    & 0.75   \\
    $n$=100                  & 0.42 & 0.43  & 0.69   & 0.66   & 0.77  & 0.74    & 0.71    & 0.83   \\
    $n$=300                  & 0.47 & 0.49  & 0.71   & 0.69   & 0.82  & 0.76    & 0.73    & 0.88   \\
    $n$=500                  & 0.48 & 0.50  & 0.73   & 0.71   & 0.83  & 0.78    & 0.76    & 0.89   \\ \hline
    \end{tabular}
\end{footnotesize}
\caption{Finite sample efficiencies of several scatter matrices: $p=2$}
\label{table:FSEtable2}
\end{table}

\begin{table}
\begin{footnotesize}
   \begin{tabular}{c|cc|ccc|ccc}
    \hline
    3-variate $t_5$    & SCM  & Tyler & HSD-CM & MhD-CM & PD-CM & HSD-wCM & MhD-wCM & PD-wCM \\ \hline
    $n$=20             & 0.96 & 0.97  & 1.06   & 1.03   & 0.99  & 1.07    & 1.06    & 0.97   \\
    $n$=50             & 1.07 & 1.08  & 1.28   & 1.20   & 1.18  & 1.33    & 1.23    & 1.20   \\
    $n$=100            & 1.12 & 1.15  & 1.49   & 1.31   & 1.40  & 1.57    & 1.38    & 1.48   \\
    $n$=300            & 1.49 & 1.54  & 2.09   & 1.82   & 2.07  & 2.19    & 1.93    & 2.18   \\
    $n$=500            & 1.60 & 1.66  & 2.18   & 1.87   & 2.21  & 2.27    & 1.95    & 2.30   \\ \hline
    3-variate $t_6$    & SCM  & Tyler & HSD-CM & MhD-CM & PD-CM & HSD-wCM & MhD-wCM & PD-wCM \\ \hline
    $n$=20             & 0.90 & 0.92  & 1.00   & 0.99   & 0.95  & 1.02    & 1.01    & 0.94   \\
    $n$=50             & 0.95 & 0.96  & 1.16   & 1.09   & 1.09  & 1.21    & 1.14    & 1.11   \\
    $n$=100            & 0.98 & 0.99  & 1.32   & 1.22   & 1.25  & 1.38    & 1.27    & 1.29   \\
    $n$=300            & 1.10 & 1.14  & 1.57   & 1.40   & 1.58  & 1.62    & 1.47    & 1.64   \\
    $n$=500            & 1.17 & 1.20  & 1.57   & 1.43   & 1.60  & 1.63    & 1.51    & 1.67   \\ \hline
    3-variate $t_{10}$ & SCM  & Tyler & HSD-CM & MhD-CM & PD-CM & HSD-wCM & MhD-wCM & PD-wCM \\ \hline
    $n$=20             & 0.87 & 0.88  & 0.95   & 0.94   & 0.90  & 0.97    & 0.98    & 0.89   \\
    $n$=50             & 0.77 & 0.79  & 0.96   & 0.92   & 0.94  & 0.99    & 0.96    & 0.95   \\
    $n$=100            & 0.75 & 0.76  & 1.02   & 0.95   & 1.01  & 1.06    & 1.00    & 1.05   \\
    $n$=300            & 0.73 & 0.75  & 1.03   & 0.98   & 1.10  & 1.08    & 1.03    & 1.15   \\
    $n$=500            & 0.73 & 0.76  & 1.02   & 0.98   & 1.09  & 1.06    & 1.02    & 1.14   \\ \hline
    3-variate $t_{15}$ & SCM  & Tyler & HSD-CM & MhD-CM & PD-CM & HSD-wCM & MhD-wCM & PD-wCM \\ \hline
    $n$=20             & 0.84 & 0.86  & 0.92   & 0.92   & 0.89  & 0.94    & 0.94    & 0.87   \\
    $n$=50             & 0.75 & 0.76  & 0.92   & 0.90   & 0.90  & 0.96    & 0.94    & 0.93   \\
    $n$=100            & 0.66 & 0.67  & 0.91   & 0.87   & 0.95  & 0.96    & 0.92    & 1.00   \\
    $n$=300            & 0.61 & 0.64  & 0.90   & 0.87   & 1.00  & 0.93    & 0.91    & 1.04   \\
    $n$=500            & 0.65 & 0.67  & 0.89   & 0.87   & 0.99  & 0.93    & 0.91    & 1.03   \\ \hline
    3-variate $t_{25}$ & SCM  & Tyler & HSD-CM & MhD-CM & PD-CM & HSD-wCM & MhD-wCM & PD-wCM \\ \hline
    $n$=20             & 0.78 & 0.79  & 0.87   & 0.89   & 0.87  & 0.89    & 0.92    & 0.86   \\
    $n$=50             & 0.70 & 0.71  & 0.88   & 0.86   & 0.88  & 0.91    & 0.90    & 0.90   \\
    $n$=100            & 0.61 & 0.63  & 0.86   & 0.83   & 0.89  & 0.90    & 0.88    & 0.94   \\
    $n$=300            & 0.58 & 0.59  & 0.83   & 0.80   & 0.92  & 0.87    & 0.85    & 0.98   \\
    $n$=500            & 0.62 & 0.64  & 0.83   & 0.82   & 0.94  & 0.88    & 0.87    & 0.99   \\ \hline
    3-variate Normal   & SCM  & Tyler & HSD-CM & MhD-CM & PD-CM & HSD-wCM & MhD-wCM & PD-wCM \\ \hline
    $n$=20             & 0.76 & 0.78  & 0.85   & 0.87   & 0.84  & 0.87    & 0.90    & 0.83   \\
    $n$=50             & 0.66 & 0.67  & 0.82   & 0.81   & 0.84  & 0.86    & 0.86    & 0.86   \\
    $n$=100            & 0.56 & 0.58  & 0.77   & 0.75   & 0.83  & 0.82    & 0.79    & 0.87   \\
    $n$=300            & 0.53 & 0.55  & 0.75   & 0.74   & 0.85  & 0.79    & 0.78    & 0.90   \\
    $n$=500            & 0.56 & 0.58  & 0.76   & 0.76   & 0.87  & 0.80    & 0.80    & 0.92   \\ \hline
    \end{tabular}
\end{footnotesize}
\caption{Finite sample efficiencies of several scatter matrices: $p=3$}
\label{table:FSEtable3}
\end{table}

\begin{table}
\begin{footnotesize}
    \begin{tabular}{c|cc|ccc|ccc}
    \hline
    4-variate $t_5$    & SCM  & Tyler & HSD-CM & MhD-CM & PD-CM & HSD-wCM & MhD-wCM & PD-wCM \\ \hline
    $n$=20             & 1.04 & 1.02  & 1.10   & 1.07   & 1.02  & 1.09    & 1.07    & 0.98   \\
    $n$=50             & 1.08 & 1.08  & 1.16   & 1.16   & 1.13  & 1.19    & 1.19    & 1.13   \\
    $n$=100            & 1.31 & 1.31  & 1.42   & 1.38   & 1.36  & 1.46    & 1.44    & 1.36   \\
    $n$=300            & 1.46 & 1.54  & 1.81   & 1.76   & 1.95  & 1.88    & 1.88    & 1.95   \\
    $n$=500            & 1.92 & 1.93  & 2.23   & 2.03   & 2.31  & 2.35    & 2.19    & 2.39   \\ \hline
    4-variate $t_6$    & SCM  & Tyler & HSD-CM & MhD-CM & PD-CM & HSD-wCM & MhD-wCM & PD-wCM \\ \hline
    $n$=20             & 1.00 & 1.05  & 1.03   & 1.05   & 1.00  & 1.04    & 1.04    & 0.95   \\
    $n$=50             & 1.03 & 1.01  & 1.13   & 1.12   & 1.11  & 1.19    & 1.17    & 1.10   \\
    $n$=100            & 1.08 & 1.12  & 1.25   & 1.23   & 1.27  & 1.24    & 1.25    & 1.22   \\
    $n$=300            & 1.34 & 1.36  & 1.64   & 1.52   & 1.60  & 1.67    & 1.61    & 1.68   \\
    $n$=500            & 1.26 & 1.34  & 1.55   & 1.49   & 1.60  & 1.65    & 1.61    & 1.69   \\ \hline
    4-variate $t_{10}$ & SCM  & Tyler & HSD-CM & MhD-CM & PD-CM & HSD-wCM & MhD-wCM & PD-wCM \\ \hline
    $n$=20             & 0.90 & 0.89  & 0.95   & 0.98   & 0.98  & 0.96    & 1.01    & 0.95   \\
    $n$=50             & 0.90 & 0.91  & 1.01   & 0.98   & 0.98  & 1.03    & 1.04    & 0.99   \\
    $n$=100            & 0.87 & 0.87  & 0.93   & 0.95   & 1.01  & 0.99    & 1.01    & 1.05   \\
    $n$=300            & 0.87 & 0.87  & 1.09   & 1.09   & 1.17  & 1.14    & 1.16    & 1.23   \\
    $n$=500            & 0.88 & 0.92  & 1.10   & 1.10   & 1.23  & 1.19    & 1.22    & 1.29   \\ \hline
    4-variate $t_{15}$ & SCM  & Tyler & HSD-CM & MhD-CM & PD-CM & HSD-wCM & MhD-wCM & PD-wCM \\ \hline
    $n$=20             & 0.92 & 0.90  & 0.94   & 0.94   & 0.96  & 0.95    & 0.97    & 0.89   \\
    $n$=50             & 0.82 & 0.83  & 0.88   & 0.91   & 0.93  & 0.88    & 0.93    & 0.93   \\
    $n$=100            & 0.84 & 0.87  & 0.92   & 0.95   & 1.00  & 0.93    & 0.96    & 1.00   \\
    $n$=300            & 0.73 & 0.75  & 0.96   & 0.99   & 1.10  & 1.00    & 1.06    & 1.12   \\
    $n$=500            & 0.73 & 0.76  & 0.95   & 0.96   & 1.06  & 0.94    & 0.97    & 1.06   \\ \hline
    4-variate $t_{25}$ & SCM  & Tyler & HSD-CM & MhD-CM & PD-CM & HSD-wCM & MhD-wCM & PD-wCM \\ \hline
    $n$=20             & 0.89 & 0.92  & 0.92   & 0.92   & 0.90  & 0.96    & 0.95    & 0.89   \\
    $n$=50             & 0.82 & 0.84  & 0.89   & 0.90   & 0.91  & 0.93    & 0.96    & 0.92   \\
    $n$=100            & 0.77 & 0.76  & 0.90   & 0.90   & 0.96  & 0.94    & 0.98    & 1.04   \\
    $n$=300            & 0.73 & 0.77  & 0.93   & 0.91   & 0.98  & 1.00    & 0.98    & 1.03   \\
    $n$=500            & 0.67 & 0.71  & 0.83   & 0.83   & 0.96  & 0.88    & 0.90    & 1.00   \\ \hline
    4-variate Normal   & SCM  & Tyler & HSD-CM & MhD-CM & PD-CM & HSD-wCM & MhD-wCM & PD-wCM \\ \hline
    $n$=20             & 0.82 & 0.84  & 0.87   & 0.90   & 0.91  & 0.89    & 0.93    & 0.89   \\
    $n$=50             & 0.80 & 0.81  & 0.87   & 0.88   & 0.88  & 0.88    & 0.92    & 0.88   \\
    $n$=100            & 0.68 & 0.71  & 0.80   & 0.85   & 0.91  & 0.82    & 0.86    & 0.92   \\
    $n$=300            & 0.61 & 0.63  & 0.82   & 0.85   & 0.93  & 0.86    & 0.91    & 0.96   \\
    $n$=500            & 0.60 & 0.64  & 0.77   & 0.80   & 0.90  & 0.82    & 0.86    & 0.96   \\ \hline
    \end{tabular}
\end{footnotesize}
\caption{Finite sample efficiencies of several scatter matrices: $p=4$}
\label{table:FSEtable4}
\end{table}

\section{Examples in real data analysis}\label{section:sec6}

\subsection{Bus data}
This dataset is available in the R package \texttt{rrcov}, and consists of data on images of 218 buses. The 18 variables here correspond to several features related to these images. Here we extend upon the analysis in \cite{maronna06}, pp. 213 to compare the classical PCA and 3 different methods of robust PCA (including that using SCM) with our DCM-based method. Similar to the original analysis, we set aside variable 9 and scale the other variables by dividing with their respective median absolute deviations (MAD). This is done because all the variables had much larger standard deviations compared to their MADs, and variable 9 had MAD = 0.

For the sake of uniformity, we use projection depth as our fixed depth function while doing depth-based PCA in our data analysis examples. We compare the outputs of classical (CPCA) and depth-based PCA (DPCA) with the following 3 robust methods: Spherical PCA, i.e. PCA based on the SCM (SPCA), PCs obtained by the ROBPCA algorithm \citep{hubert05}, and the Minimum Covariance Determinant (MCD) estimator \citep{Rousseeuw84leastmedian} (MPCA). Table \ref{table:bus_table1} gives the proportions of variability that are left unexplained after the top $q$ $(= 1,...,6)$ components are taken into account in each of the 5 methods. The first PC in classical PCA seems to explain a much higher proportion of variability in original data than robust methods. However, as noted in \cite{maronna06} and \cite{hubert05}, this is a result of the classical variances being inflated due to outliers in the direction of the first principal axis. Among the robust methods, the proportions of unexplained variances are highest for DCM-based PCA for all values of $q$.

Table \ref{table:bus_table2} demonstrates why robust methods actually give a better representation of the underlying data structure than classical PCA here. Each of its column lists different quantiles of the squared orthogonal distance for a sample point from the hyperplane formed by top 3 PCs estimated by the corresponding method. For PCA based on projection-DCM, the estimated principal component subspaces are closer to the data than CPCA for more than 90\% of samples, and the distance only becomes larger for higher quantiles. This means that for CPCA, estimated basis vectors of the hyperspace get pulled by extreme outlying points, while the influence of these outliers is very low for DPCA. SPCA and ROBPCA perform very closely in this respect, the percentage of points that have less squared distance than CPCA being between 80\% and 90\% for both of them. This percentage is only 50\% for MPCA, which suggests that the corresponding 3-dimensional subspace estimated by MCD is possibly not an accurate representation of the truth.

\begin{table}[t]
\centering
    \begin{tabular}{c|ccccc}
    \hline
    \begin{large} $q$ \end{large} & \multicolumn{5}{c}{Method of PCA}         \\ \cline{2-6}
    ~                   & CPCA     & SPCA & ROBPCA & MPCA & DPCA \\\hline 
    1                   & 0.188         & 0.549     & 0.410  & 0.514     & 0.662     \\
    2                   & 0.084         & 0.272     & 0.214  & 0.337     & 0.359     \\
    3                   & 0.044         & 0.182     & 0.121  & 0.227     & 0.237     \\
    4                   & 0.026         & 0.135     & 0.083  & 0.154     & 0.173     \\
    5                   & 0.018         & 0.099     & 0.054  & 0.098     & 0.115     \\
    6                   & 0.012         & 0.069     & 0.036  & 0.070     & 0.084     \\ \hline
    \end{tabular}
    \caption{Unexplained proportions of variability by PCA models with $q$ components for bus data}
    \label{table:bus_table1}
\end{table}

\begin{table}[t]
\centering
    \begin{tabular}{c|ccccc}
    \hline
    Quantile & \multicolumn{5}{c}{Method of PCA}         \\ \cline{2-6}
    ~                   & CPCA     & SPCA & ROBPCA & MPCA & DPCA \\\hline 
    10\%      & 1.9       & 1.2       & 1.2    & 1.0       & 1.2       \\
    20\%      & 2.3       & 1.6       & 1.6    & 1.3       & 1.6       \\
    30\%      & 2.8       & 1.8       & 1.8    & 1.7       & 1.9       \\
    40\%      & 3.2       & 2.2       & 2.1    & 2.1       & 2.3       \\
    50\%      & 3.7       & 2.6       & 2.5    & 3.1       & 2.6       \\
    60\%      & 4.4       & 3.1       & 3.0    & 5.9       & 3.2       \\
    70\%      & 5.4       & 3.8       & 3.9    & 25.1      & 3.9       \\
    80\%      & 6.5       & 5.2       & 4.8    & 86.1      & 4.8       \\
    90\%      & 8.2       & 9.0       & 10.9   & 298.2     & 6.9      \\
    Max       & 24        & 1037      & 1055   & 1037      & 980      \\\hline
    \end{tabular}
    \caption{Quantiles to squared distance from 3-principal component hyperplanes for bus data}
    \label{table:bus_table2}
\end{table}

\subsection{Octane data}
\begin{figure}[t]
	\centering
		\includegraphics[height=6cm]{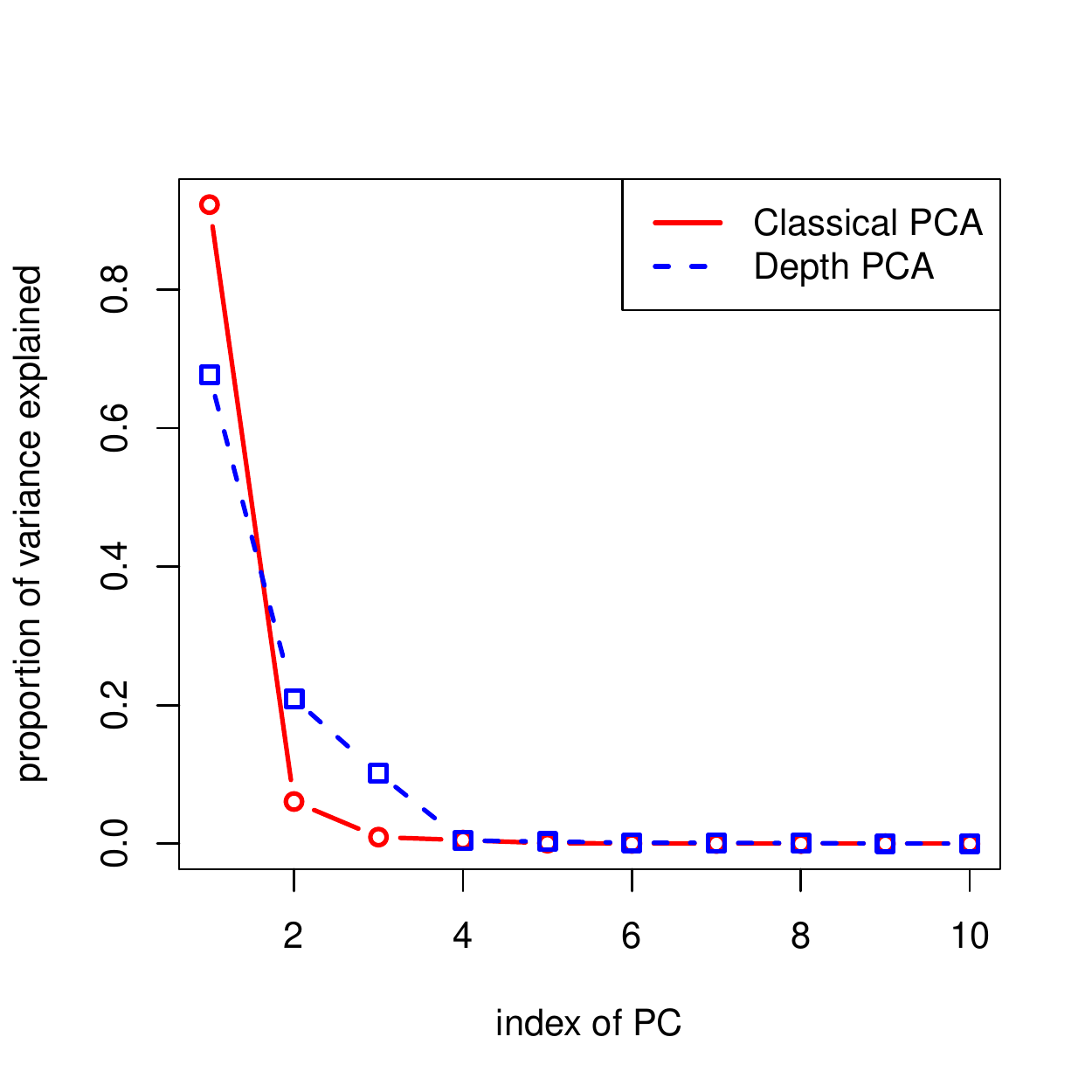}
	\caption{Explained variance proportions for two types of PCA on octane data}
	\label{fig:screeplot}
\end{figure}

We now apply our method to a high-dimensional dataset and demonstrate its effectiveness in outlier detection. Due to \cite{esbensen94}, this dataset consists of 226 variables and 39 observations. Each observation is a gasoline sample with a certain octane number, and have their NIR absorbance spectra measured in 2 nm intervals between 1100 - 1550 nm. There are 6 outliers here: compounds 25, 26 and 36-39, which contain alcohol. Proportions of \textit{explained} variability by PCs obtained by the two methods are given in Figure \ref{fig:screeplot}. Once again, the first PC in CPCA explains a much larger proportion of variance than DPCA. Second and third PCs obtained by DPCA explain higher proportions of variability than the corresponding components in CPCA.

\begin{figure}[t]
	\centering
		\includegraphics[height=6.5cm]{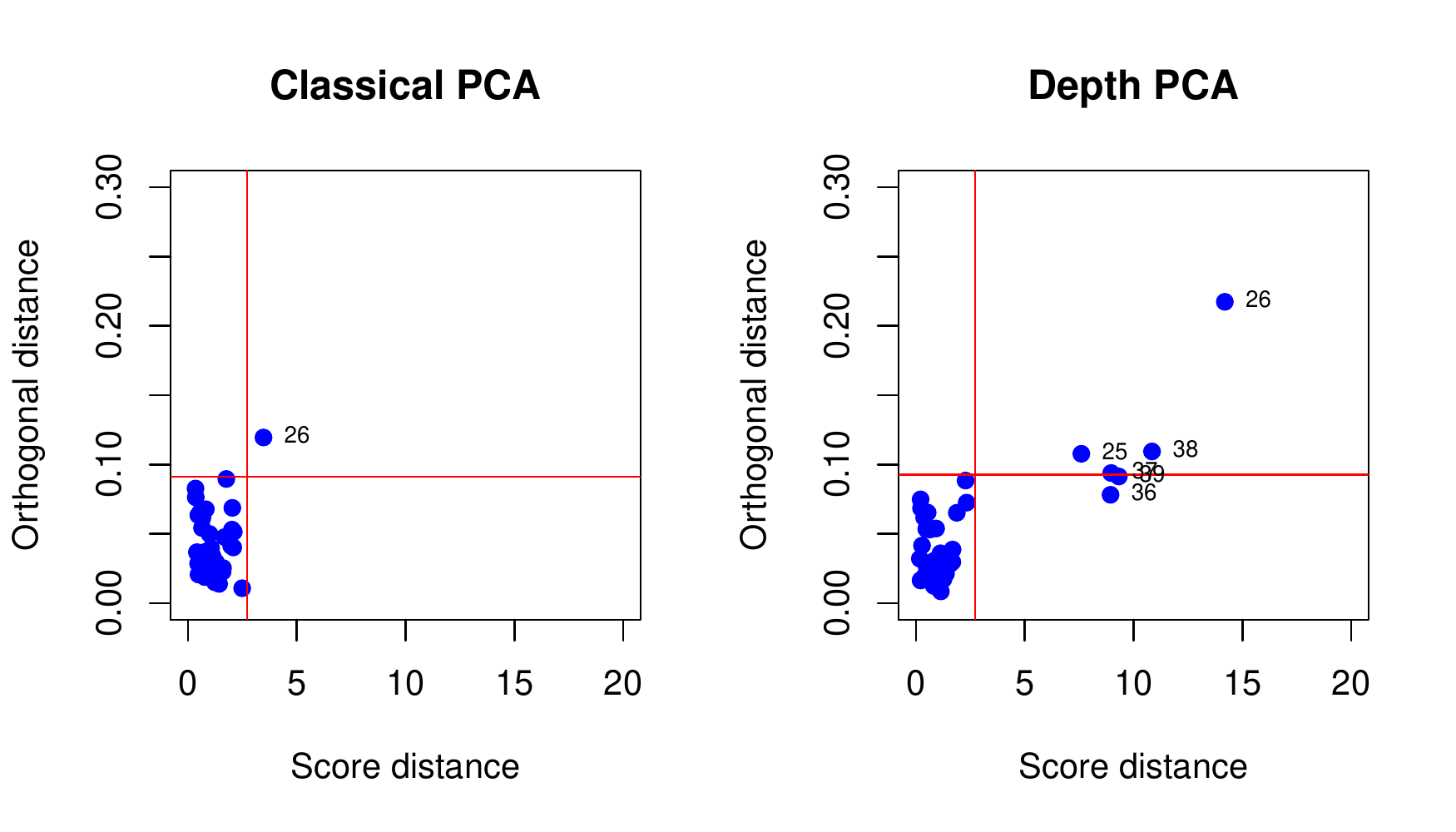}
	\caption{Distance plots for two types of PCA on octane data}
	\label{fig:distplot}
\end{figure}

In both methods, the first two PCs explains a large amount (98\% for CPCA, 89\% for DPCA) of underlying variability. But in DPCA, these PCs are more effective in detecting outliers, which we demonstrate in Figure \ref{fig:distplot}. For any method of PCA with $k$ components on a dataset of $n$ observations and $p$ variables, the \textbf{score distance} (SD) and \textbf{orthogonal distance} (OD) for $i^\text{th}$ observation ($i=1,2,...,n$) are defined as:
$$ SD_i = \sqrt{ \sum_{j=1}^k \frac{s^2_{ij}}{\lambda_j}}; \quad OD_i = \| \bfx_i - P\bfs_i^T \| $$
where $S = (\bfs_1, \ldots ,\bfs_n)^T$ is the $n\times k$ scoring matrix, $P$ the $p\times k$ loading matrix, and $\lambda_1,\ldots ,\lambda_k$ are eigenvalues obtained from the PCA, and $\bfx_1,\ldots,\bfx_n$ are the $n$ observation vectors. We only consider the first 2 PCs here, so $k$ is set to 2. From a practical standpoint, $SD_i$ can be interpreted as a weighted norm of the projection of the $i^\text{th}$ point on the hyperplane formed by first $k$ principal components, and $OD_i$ the orthogonal distance of point $i$ from that hyperplane. For outlier detection, following \cite{hubert05} we set the upper cutoff values for score distances at $\sqrt{\chi^2_{2,.975}}$ and orthogonal distances at $[\text{median}(OD^{2/3}) + \text{MAD}(OD^{2/3})\Phi^{-1}(0.975)]^{3/2}$, where $\Phi(.)$ is the standard normal cumulative distribution function. These cutoffs are marked by red lines in the diagnostic plots in figure \ref{fig:distplot}. In the figure we see that CPCA detects only 1 of the 6 outliers (point 26, left panel). Compared to this the diagnostic plot corresponding to DPCA in the right panel correctly detects all 6 outlying points. All of them have larger score distances than the score cutoff value, while points 25, 26 and 38 have higher-than-cutoff orthogonal distances as well. Interestingly, the arrangement of outlying points here is similar to that in the diagnostic plot corresponding to ROBPCA, which can be found in figure 5(b) of \cite{hubert05}.


\section{Conclusion}\label{section:sec7}

In the above sections we introduce a covariance matrix based on depth-based multivariate ranks that keeps the eigenvectors of the actual population unchanged for elliptical distributions. We provide asymptotic results for the sample DCM, its eigenvalues and eigenvectors. Bounded influence functions as well as simulation studies suggest that the eigenvector estimates obtained from the DCM are highly robust, yet do not lose much in terms of efficiency. Thus it provides a plausible alternative to existing approaches of robust PCA that are based on estimation of covariance matrices (for example SCM, Tyler's scatter matrix, D\"{u}mbgen's symmetrized shape matrix).

An immediate extension of this would be to study the depth-weighted iterated scatter matrices, i.e. matrices $\Sigma_{Dw}$ that are solution to the following type of equations:

$$ \Sigma_{Dw} = E \left[ \frac{(\tilde D_\bfX(\bfx))^2 (\bfx - \bfmu) (\bfx - \bfmu)^T}{(\bfx - \bfmu)^T \Sigma_{Dw}^{-1} (\bfx - \bfmu)} \right] $$
as their eigenvector estimates seem to have better efficiency than those obtained from the corresponding DCMs. Unlike the DCM, these matrices will possess the affine equivariance property. It is possible to develop tests for central and elliptic symmetry based on the decomposition of the multivariate rank vector in equation \ref{equation:rankdecomp}. The result in Theorem \ref{Theorem:covform} is based on the fact that $E(\tilde \bfX) = \bf0$, and it holds for any centrally symmetric underlying distribution. Moreover, the depth of a point in the standardized scale (i.e. $\bfz$-scale) does not depend on the direction of $\bfz$. This is not possible without circular symmetry of $\bfz$, so any test of independence between $D_\bfZ (\bfz)$ and $\bfS (\bfz)$ can be seen as a test of ellipticity for the original random variable $\bfX$. Finally the applicability of this procedure in high-dimensional and functional data remains to be explored.

\appendix
\section*{Appendix}
\numberwithin{equation}{section}
\section{\textbf{Form of $V_{D,S}(F)$}}\label{section:appA}

First observe that for $F$ having covariance matrix $\Sigma = \Gamma\Lambda\Gamma^T$,
$$ V_{D,S}(F)  = (\Gamma \otimes \Gamma) V_{D,S}(F_\Lambda) (\Gamma \otimes \Gamma)^T$$
where $F_\Lambda$ has the same elliptic distribution as $F$, but with covariance matrix $\Lambda$. Now,
\begin{eqnarray*}
V_{D,S} (F_\Lambda) &=& E \left[ vec \left\{ \frac{(\tilde D_\bfZ (\bfz))^2 \Lambda^{1/2} \bfz\bfz^T \Lambda^{1/2}}{\bfz^T\Lambda\bfz} - \Lambda_{D,S} \right\} vec^T \left\{ \frac{(\tilde D_\bfZ (\bfz))^2 \Lambda^{1/2} \bfz\bfz^T \Lambda^{1/2}}{\bfz^T\Lambda\bfz} - \Lambda_{D,S} \right\} \right]\\
&=& E \left[ vec \left\{ (\tilde D_\bfZ (\bfz))^2 SS(\Lambda^{1/2}\bfz; \bf0) \right\} vec^T \left\{ (\tilde D_\bfZ (\bfz))^2 SS(\Lambda^{1/2}\bfz; \bf0) \right\} \right]\\
&& - \quad vec(\Lambda_{D,S}) vec^T(\Lambda_{D,S})
\end{eqnarray*}

The matrix $vec(\Lambda_{D,S}) vec^T(\Lambda_{D,S})$ consists of elements $\lambda_i\lambda_j$ at $(i,j)^\text{th}$ position of the $(i,j)^\text{th}$ block, and 0 otherwise. These positions correspond to variance and covariance components of on-diagonal elements. For the expectation matrix, all its elements are of the form $E[\sqrt{\lambda_a \lambda_b \lambda_c \lambda_d} z_a z_b z_c z_d . (\tilde D_\bfZ (\bfz))^4 / (\bfz^T \Lambda \bfz)^2]$, with $1 \leq a,b,c,d \leq p$. Since $(\tilde D_\bfZ (\bfz))^4 / (\bfz^T \Lambda \bfz)^2$ is even in $\bfz$, which has a circularly symmetric distribution, all such expectations will be 0 unless $a=b=c=d$, or they are pairwise equal. Following a similar derivation for spatial sign covariance matrices in \cite{magyar14}, we collect the non-zero elements and write the matrix of expectations:
$$ (I_{p^2} + K_{p,p}) \left\{ \sum_{a=1}^p \sum_{b=1}^p \gamma^D_{ab} (\bfe_a \bfe_a^T \otimes  \bfe_b \bfe_b^T) - \sum_{a=1}^p \gamma^D_{aa} (\bfe_a \bfe_a^T \otimes  \bfe_a \bfe_a^T) \right\} + \sum_{a=1}^p \sum_{b=1}^p \gamma^D_{ab} (\bfe_a \bfe_b^T \otimes  \bfe_a \bfe_b^T) $$
where $I_k = (\bfe_1,...,\bfe_k), K_{m,n} = \sum_{i=1}^m \sum_{j=1}^n J_{ij} \otimes J_{ij}^T$ with $J_{ij}$ the $m \times n$ matrix having 1 as $(i,j)^\text{th}$ element and 0 elsewhere, and $\gamma^D_{mn} = E[ \lambda_m \lambda_n z_m^2 z_n ^2 . (\tilde D_\bfZ (\bfz))^4 / (\bfz^T \Lambda \bfz)^2]; 1 \leq m,n \leq p$.

\paragraph{}Putting everything together, denote $\hat S^D(F_\Lambda) = \sum_{i=1}^n (\tilde D^n_\bfZ (\bfz_i))^2 SS(\Lambda^{1/2}\bfz_i; \hat \bfmu_n)/n $. Then the different types of elements in the matrix $V_{D,S}(F_\Lambda)$ are as given below ($1 \leq a,b,c,d \leq p$):

\begin{itemize}
\item Variance of on-diagonal elements
$$ AVar( \sqrt n \hat S^D_{aa} (F_\Lambda)) = E \left[ \frac{(\tilde D_\bfZ (\bfz))^4 \lambda_a^2 z_a^4}{(\bfz^T \Lambda \bfz)^2} \right] - \lambda_{D,S,a}^2 $$

\item Variance of off-diagonal elements ($a \neq b$)
$$ AVar( \sqrt n \hat S^D_{ab} (F_\Lambda)) = E \left[ \frac{(\tilde D_\bfZ (\bfz))^4 \lambda_a \lambda_b z_a^2 z_b^2}{(\bfz^T \Lambda \bfz)^2} \right] $$

\item Covariance of two on-diagonal elements ($a \neq b$)
$$ ACov(\sqrt n \hat S^D_{aa} (F_\Lambda), \sqrt n \hat S^D_{bb} (F_\Lambda))
= E \left[ \frac{(\tilde D_\bfZ (\bfz))^4 \lambda_a \lambda_b z_a^2 z_b^2}{(\bfz^T \Lambda \bfz)^2} \right] - \lambda_{D,S,a} \lambda_{D,S,b} $$

\item Covariance of two off-diagonal elements ($a \neq b \neq c \neq d$)
$$ ACov(\sqrt n \hat S^D_{ab} (F_\Lambda), \sqrt n \hat S^D_{cd} (F_\Lambda)) = 0 $$

\item Covariance of one off-diagonal and one on-diagonal element ($a \neq b \neq c$)
$$ ACov(\sqrt n \hat S^D_{ab} (F_\Lambda), \sqrt n \hat S^D_{cc} (F_\Lambda)) = 0 $$
\end{itemize}

\section{Proofs}\label{section:appB}

\begin{proof}[Proof of Theorem  \ref{Theorem:covform}]
The proof follows directly from writing out the expression of $Cov ( \tilde \bfX)$:
\begin{eqnarray*}
Cov(\tilde\bfX) &=& E(\tilde\bfX \tilde\bfX^T) - E(\tilde\bfX) E(\tilde\bfX)^T\\
&=& \Gamma . E \left[ (\tilde D_\bfZ(\bfz))^2 \frac{\|\bfz\|^2}{\| \Lambda^{1/2}\bfz\|} \Lambda^{1/2} \bfS(\bfz) \bfS(\bfz)^T \Lambda^{1/2} \right] \Gamma^T - {\bf 0}_p {\bf 0}_p^T\\
&=& \Gamma .E \left[ (\tilde D_\bfZ(\bfz))^2 \frac{\Lambda^{1/2} \bfz \bfz^T \Lambda^{1/2}}{\bfz^T \Lambda \bfz} \right] \Gamma^T
\end{eqnarray*}
\end{proof}

\begin{proof}[Proof of Lemma \ref{Lemma:lemma1}]
For two positive definite matrices $A,B$, we denote by $A>B$ that $A-B$ is positive definite. Also, denote
$$ S_n = \sqrt n \left[ \frac{1}{n} \sum_{i=1}^n \left| (\tilde D^{n} _\bfX (\bfx_i))^2  - (\tilde D_\bfX (\bfx_i))^2 \right| SS(\bfx_i; \hat\bfmu_n) \right] $$
%
Now due to the assumption of uniform convergence, given $\epsilon>0$ we can find $N \in \mathbb{N}$ such that
\begin{equation}
\label{equation:lemma1eq}
| (\tilde D^{n_1}_\bfX(\bfx_i))^2 - (\tilde D_\bfX(\bfx_i))^2 | < \epsilon
\end{equation}
for all $n_1 \geq N; i = 1,2,...,n_1$. This implies
\begin{eqnarray}
\label{eqn:lemma1eq2}
S_{n_1} &< & \epsilon \sqrt{n_1} \left[ \frac{1}{n_1} \sum_{i=1}^{n_1} SS(\bfx_i; \hat\bfmu_{n_1}) \right]\notag\\
&=& \epsilon \sqrt{n_1} \left[ \frac{1}{n_1} \sum_{i=1}^{n_1} \left\{ SS(\bfx_i; \hat\bfmu_{n_1}) - SS(\bfx_i; \bfmu) \right\} + \frac{1}{n_1} \sum_{i=1}^{n_1} SS(\bfx_i; \bfmu) \right]
\end{eqnarray}

We now construct a sequence of positive definite matrices $\{A_k (B_k+C_k) : k \in \mathbb N\} $ so that
$$ A_k = \frac{1}{k}, \quad B_k = \sqrt{N_k} \left[ \frac{1}{N_k} \sum_{i=1}^{N_k} \left\{ SS(\bfx_i; \hat\bfmu_{N_k}) - SS(\bfx_i; \bfmu) \right\} \right]$$
$$\quad C_k = \sqrt{N_k} \left[ \frac{1}{N_k} \sum_{i=1}^{N_k} SS(\bfx_i; \bfmu) \right] $$
where $N_k \in \mathbb N$ gives the relation (\ref{equation:lemma1eq}) in place of $N$ when we take $\epsilon = 1/k$. Under conditions $ E\|\bfx - \bfmu\|^{-3/2} < \infty $ and $\sqrt n (\hat\bfmu_n - \bfmu) = O_P(1)$, the sample SCM with unknown location parameter $\hat\bfmu_n$ has the same asymptotic distribution as the SCM with known location $\bfmu$ \citep{durre14}, hence $B_k = o_P(1)$, thus $A_k (B_k+C_k) \stackrel{P}{\rightarrow} 0$.

Now (\ref{eqn:lemma1eq2}) implies that for any $\epsilon_1 > 0$, $S_{N_k} > \epsilon_1 \Rightarrow A_k (B_k + C_k) > \epsilon_1$, which means $ P(S_{N_k} > \epsilon_1) < P(A_k (B_k + C_k) > \epsilon_1)$. Hence the subsequence $\{S_{N_k}\} \stackrel{P}{\rightarrow} 0$. Since the main sequence $\{S_k\}$ is bounded below by 0, this implies $\{S_k\} \stackrel{P}{\rightarrow} 0$. Finally, we have that
\begin{eqnarray}
\sqrt n \left[
\frac{1}{n} \sum_{i=1}^n (\tilde D^n_\bfX (\bfx_i))^2 SS(\bfx_i; \hat\bfmu_n) -
\frac{1}{n} \sum_{i=1}^n (\tilde D_\bfX (\bfx_i))^2 SS(\bfx_i; \bfmu) \right] &\leq & \hspace{1em} \notag\\
S_n +  \sqrt{n} \left[ \frac{1}{n} \sum_{i=1}^{n} \left\{ SS(\bfx_i; \hat\bfmu_{n}) - SS(\bfx_i; \bfmu) \right\} \right] &&
\end{eqnarray}
Since the second summand on the right hand side is $o_P(1)$ due to \cite{durre14} as mentioned before, we have the needed.
\end{proof}

\begin{proof}[Proof of Theorem \ref{Theorem:rootn}]
The quantity in the statement of the theorem can be broken down as:
\begin{eqnarray*}
\sqrt n \left[ vec\left\{ \frac{1}{n} \sum_{i=1}^n (\tilde D^n_\bfX (\bfx_i))^2 SS(\bfx_i; \hat\bfmu_n) \right\} - vec\left\{ \frac{1}{n} \sum_{i=1}^n (\tilde D_\bfX (\bfx_i))^2 SS(\bfx_i; \bfmu) \right\} \right] +\\
\sqrt n \left[ vec\left\{ \frac{1}{n} \sum_{i=1}^n (\tilde D_\bfX (\bfx_i))^2 SS(\bfx_i; \bfmu) \right\} - E \left[ vec\left\{ (\tilde D_\bfX (\bfx))^2 SS(\bfx; \bfmu) \right\} \right] \right]
\end{eqnarray*}
The first part goes to 0 in probability by Lemma \ref{Lemma:lemma1}, and applying Slutsky's theorem we get the required convergence.
\end{proof}

\begin{proof}[Proof of Theorem \ref{Theorem:decomp}]
See \cite{taskinen12}
\end{proof}

\begin{proof}[Proof of Corollary \ref{Corollary:eigendist}]
In spirit, this corollary is similar to Theorem 13.5.1 in \cite{anderson}, and indeed, \cite{taskinen12} used this theorem to prove Theorem \ref{Theorem:decomp}. Due to the decomposition (\ref{equation:decompEq}) we have, for the distribution $F_\Lambda$, the following relation between any off-diagonal element of $\hat S^D(F_\Lambda)$ and the corresponding element in the estimate of eigenvectors $\hat\Gamma_D (F_\Lambda)$:

$$ \sqrt n \hat\gamma_{D,ij} (F_\Lambda) = \sqrt n \frac{\hat S^D_{ij} (F_\Lambda)}{\lambda_{D,S,i} - \lambda_{D,S,j}}; \quad i \neq j$$

So that for eigenvector estimates of the original $F$ we have

\begin{equation} \label{equation:app1}
\sqrt n (\hat\bfgamma_{D,i} - \bfgamma_i) = \sqrt n \Gamma (\hat \bfgamma_{D,i}(F_\Lambda) - \bfe_i ) = \sqrt n \left[ \sum_{k=1; k \neq i}^p \hat \gamma_{D,ik}(F_\Lambda)\bfgamma_k + (\hat \gamma_{D,ii}(F_\Lambda) - 1)\bfgamma_i \right]
\end{equation}

$\sqrt n (\hat \gamma_{D,ii}(F_\Lambda) - 1) =  o_P(1)$ and $ACov(\sqrt n \hat S^D_{ik}(F_\Lambda), \sqrt n \hat S^D_{il}(F_\Lambda)) = 0$ for $k \neq l$, so the above equation implies

$$ AVar(\bfg_i) = AVar (\sqrt n (\hat\bfgamma_{D,i} - \bfgamma_i)) = \sum_{k=1; k \neq i}^p \frac{AVar(\sqrt n \hat S^D_{ik}(F_\Lambda))}{(\lambda_{D,s,i} - \lambda_{D,S,k})^2} \bfgamma_k \bfgamma_k^T $$

For the covariance terms, from (\ref{equation:app1}) we get, for $i \neq j$,

\begin{eqnarray*}
ACov(\bfg_i, \bfg_j) &=& ACov (\sqrt n (\hat\bfgamma_{D,i} - \bfgamma_i), \sqrt n (\hat\bfgamma_{D,j} - \bfgamma_j))\\
&=& ACov \left( \sum_{k=1; k \neq i}^p \sqrt n \hat \gamma_{D,ik}(F_\Lambda)\bfgamma_k, \sum_{k=1; k \neq j}^p \sqrt n \hat \gamma_{D,jk}(F_\Lambda)\bfgamma_k \right)\\
&=& ACov \left( \sqrt n \hat \gamma_{D,ij}(F_\Lambda)\bfgamma_j, \sqrt n \hat \gamma_{D,ji}(F_\Lambda)\bfgamma_i \right)\\
&=& - \frac{AVar(\sqrt n \hat S^D_{ij}(\Lambda))}{(\lambda_{D,s,i} - \lambda_{D,S,j})^2} \bfgamma_j \bfgamma_i^T
\end{eqnarray*}

The exact forms given in the statement of the corollary now follows from the  Form of $V_{D,S}$ in Appendix \ref{section:appA}.

\paragraph{}For the on-diagonal elements of $\hat S^D(F_\Lambda)$ Theorem \ref{Theorem:decomp} gives us $ \sqrt n \hat\lambda_{D,s,i} (F_\Lambda) = \sqrt n \hat S^D_{ii}(F_\Lambda)$ for $i = 1,...,p$. Hence

\begin{eqnarray*}
AVar(l_i) &=& AVar(\sqrt n \hat\lambda_{D,s,i} - \sqrt n \lambda_{D,S,i})\\
&=& AVar(\sqrt n \hat\lambda_{D,s,i} (F_\Lambda) - \sqrt n \lambda_{D,S,i}(F_\Lambda))\\
&=& AVar(\sqrt n S^D_{ii}(F_\Lambda))
\end{eqnarray*}

A similar derivation gives the expression for $AVar(l_i,l_j); i \neq j$. Finally, since the asymptotic covariance between an on-diagonal and an off-diagonal element of $\hat S^D(F_\Lambda)$, it follows that the elements of $G$ and diagonal elements of $L$ are independent.
\end{proof}

\bibliographystyle{plainnat}
\bibliography{Depth-scatter-summary}

\begin{thebibliography}{42}
\providecommand{\natexlab}[1]{#1}
\providecommand{\url}[1]{\texttt{#1}}
\expandafter\ifx\csname urlstyle\endcsname\relax
  \providecommand{\doi}[1]{doi: #1}\else
  \providecommand{\doi}{doi: \begingroup \urlstyle{rm}\Url}\fi

\bibitem[Anderson(3rd ed. 2003)]{anderson}
T.W. Anderson.
\newblock \emph{An {I}ntroduction to {M}ultivariate {S}tatistical {A}nalysis}.
\newblock Wiley, Hoboken, NJ, 3rd ed. 2003.

\bibitem[Brown(1983)]{brown83}
B.M. Brown.
\newblock Statistical {U}se of the {S}patial {M}edian.
\newblock \emph{J. Royal Statist. Soc. B}, 45:\penalty0 25--30, 1983.

\bibitem[Chernozhukov et~al.(2014)Chernozhukov, Galichon, Hallin, and
  Henry]{Chernozhukov14}
V.~Chernozhukov, A.~Galichon, M.~Hallin, and M.~Henry.
\newblock {Monge-Kantorovich Depths, Quantiles, Ranks and Signs}.
\newblock \url{http://arxiv.org/abs/1412.8434}, 2014.

\bibitem[Croux and Haesbroeck(2000)]{croux00}
C.~Croux and G.~Haesbroeck.
\newblock {Principal Component Analysis based on Robust Estimators of the
  Covariance or Correlation Matrix: Influence Functions and Efficiencies}.
\newblock \emph{Biometrika}, 87:\penalty0 603--618, 2000.

\bibitem[D\"{u}mbgen(1992)]{Dumbgen92}
L.~D\"{u}mbgen.
\newblock Limit theorems for the simplicial depth.
\newblock \emph{Statist. Probab. Lett.}, 14:\penalty0 119--128, 1992.

\bibitem[D\"{u}mbgen(1998)]{dumbgen98}
L.~D\"{u}mbgen.
\newblock On {T}yler's {M}-functional of scatter in high dimension.
\newblock \emph{Ann. Inst. Statist. Math.}, 50:\penalty0 471--491, 1998.

\bibitem[D\"{u}rre et~al.(2014)D\"{u}rre, Vogel, and Tyler]{durre14}
A.~D\"{u}rre, D.~Vogel, and D.E. Tyler.
\newblock The spatial sign covariance matrix with unknown location.
\newblock \emph{J. Mult. Anal.}, 130:\penalty0 107--117, 2014.

\bibitem[Dutta and Ghosh(2012)]{dutta12}
S.~Dutta and A.K. Ghosh.
\newblock On robust classification using projection depth.
\newblock \emph{Ann. Inst. Stat. Math.}, 64-3:\penalty0 657--676, 2012.

\bibitem[Esbensen et~al.(1994)Esbensen, Sch\"{o}nkopf, and
  Midtgaard]{esbensen94}
K.~H. Esbensen, S.~Sch\"{o}nkopf, and T.~Midtgaard.
\newblock \emph{Multivariate {A}nalysis in {P}ractice}.
\newblock CAMO, Trondheim, Germany, 1994.

\bibitem[Ghosh and Chaudhuri(2005)]{ghosh05}
A.K. Ghosh and P.~Chaudhuri.
\newblock On {M}aximum {D}epth and {R}elated {C}lassifiers.
\newblock \emph{Scand. J. Statist.}, 32:\penalty0 327--350, 2005.

\bibitem[Haldane(1948)]{haldane48}
J.B.S. Haldane.
\newblock Note on the {M}edian of a {M}ultivariate {D}istribution.
\newblock \emph{Biometrika}, 35:\penalty0 414--415, 1948.

\bibitem[Hallin and Paindaveine(2002)]{HallinPaindaveine02}
M.~Hallin and D.~Paindaveine.
\newblock {Optimal tests for multivariate loca-tion based on interdirections
  and pseudo-Mahalanobis ranks}.
\newblock \emph{Ann. Statist.}, 30:\penalty0 1103--1133, 2002.

\bibitem[Hampel et~al.(1986)Hampel, Ronchetti, Rousseeuw, and Staehl]{hampel}
F.R. Hampel, E.M. Ronchetti, P.J. Rousseeuw, and W.A. Staehl.
\newblock \emph{Robust {S}tatistics: {T}he {A}pproach {B}ased on {I}nfluence
  {F}unctions}.
\newblock Wiley, New York, NY, 1986.

\bibitem[Huber(1977)]{huber77}
P.J. Huber.
\newblock \emph{Robust {S}tatistical {P}rcedures}.
\newblock Society for Industrial and Applied Mathematics, Philadelphia, 1977.

\bibitem[Hubert et~al.(2005)Hubert, Rousseeuw, and Branden]{hubert05}
M.~Hubert, P.~J. Rousseeuw, and K.~V. Branden.
\newblock {ROBPCA}: {A} {N}ew {A}pproach to {R}obust {P}rincipal {C}omponent
  {A}nalysis.
\newblock \emph{Technometrics}, 47-1:\penalty0 64--79, 2005.

\bibitem[Jornsten(2004)]{jornsten04}
R.~Jornsten.
\newblock Clustering and classification based on the $l_1$ depth.
\newblock \emph{J. Mult. Anal.}, 90-1:\penalty0 67--89, 2004.

\bibitem[Liu(1990)]{liu90}
R.Y. Liu.
\newblock On a notion of data depth based on random simplices.
\newblock \emph{Ann. Statist.}, 18:\penalty0 405--414, 1990.

\bibitem[Liu et~al.(1999)Liu, Parelius, and Singh]{LiuPareliusSingh99}
R.Y. Liu, J.M. Parelius, and K.~Singh.
\newblock Multivariate analysis by data depth: {D}escriptive statistics,
  graphics and inference (with discussion).
\newblock \emph{Ann. Statist.}, 27:\penalty0 783--858, 1999.

\bibitem[Locantore et~al.(1999)Locantore, Marron, Simpson, Tripoli, Zhang, and
  Cohen]{locantore99}
N.~Locantore, J.S. Marron, D.G. Simpson, N.~Tripoli, J.T. Zhang, and K.L.
  Cohen.
\newblock Robust principal components of functional data.
\newblock \emph{TEST}, 8:\penalty0 1--73, 1999.

\bibitem[Magyar and Tyler(2014)]{magyar14}
A.F. Magyar and D.E. Tyler.
\newblock The asymptotic inadmissibility of the spatial sign covariance matrix
  for elliptically symmetric distributions.
\newblock \emph{Biometrika}, 101:\penalty0 673--688, 2014.

\bibitem[Maronna et~al.(2006)Maronna, Martin, and Yohai]{maronna06}
R.~Maronna, D.~Martin, and V.~Y. Yohai.
\newblock \emph{Robust {S}tatistics: {T}heory and {M}ethods}.
\newblock Wiley, New York, NY, 2006.

\bibitem[Maronna et~al.(1976)Maronna, Staehl, and Yohai]{maronna76}
R.A. Maronna, W.A. Staehl, and V.Y. Yohai.
\newblock Bias-{R}obust {E}stimators of {M}ultivariate {S}catter {B}ased on
  {P}rojections.
\newblock \emph{J. Mult. Anal.}, 42:\penalty0 141--161, 1976.

\bibitem[Miao and Ben-Israel(1992)]{miao92}
J.M. Miao and A.~Ben-Israel.
\newblock On principal angles between subspaces in $\mathbb{R}^n$.
\newblock \emph{Lin. Algeb. Applic.}, 171:\penalty0 81--98, 1992.

\bibitem[M\"{o}tt\"{o}nen and Oja(1995)]{MottonenOja95}
J.~M\"{o}tt\"{o}nen and H.~Oja.
\newblock Multivariate spatial sign and rank methods.
\newblock \emph{J. Nonparametric Stat.}, 1995.

\bibitem[Oja(1983)]{oja83}
H.~Oja.
\newblock Descriptive {S}tatistics for {M}ultivariate {D}istributions.
\newblock \emph{Statist. and Prob. Lett.}, 1:\penalty0 327--332, 1983.

\bibitem[Ollilia et~al.(2003)Ollilia, Oja, and Croux]{ollilia03}
E.~Ollilia, H.~Oja, and C.~Croux.
\newblock The affine equivariant sign covariance matrix: asymptotic behavior
  and efficiencies.
\newblock \emph{J. Mult. Anal.}, 87:\penalty0 328--355, 2003.

\bibitem[Puri and Sen(1971)]{PuriSenBook}
M.~L. Puri and P.~K. Sen.
\newblock \emph{{Nonparametric Methods in Multivariate Analysis}}.
\newblock Wiley, New York, NY, 1971.

\bibitem[Rainer and Mozharovskyi(2014)]{rainerArxiv}
D.~Rainer and P.~Mozharovskyi.
\newblock Exact computation of the halfspace depth.
\newblock \url{http://arxiv.org/abs/1411.6927}, 2014.

\bibitem[Rousseeuw(1985)]{rousseeuw85}
P.~Rousseeuw.
\newblock Multivariate estimation with high breakdown point.
\newblock In \emph{Mathematical Statistics and Applications, Volume B (Proc.
  4th Pannonian Symp. Math. Statist., Bad Tatzmannsdorf, 1983)}, pages 283--97,
  Dordrecht, 1985. D. Reidel.

\bibitem[Rousseeuw(1984)]{Rousseeuw84leastmedian}
Peter~J. Rousseeuw.
\newblock Least {M}edian of {S}quares {R}egression.
\newblock \emph{J. Amer. Statist. Assoc.}, 79:\penalty0 871--880, 1984.

\bibitem[Serfling(2006)]{serfling2006}
R.~Serfling.
\newblock Depth {F}unctions in {N}onparametric {M}ultivariate {I}nference.
\newblock In \emph{DIMACS {S}eries in {D}iscrete {M}athematics and
  {T}heoretical {C}omputer {S}cience}, volume~72, pages 1--16, 2006.

\bibitem[Sguera et~al.(2014)Sguera, Galeano, and Lillo]{sguera14}
C.~Sguera, P.~Galeano, and R.~Lillo.
\newblock Spatial depth-based classification for functional data.
\newblock \emph{TEST}, 23-4:\penalty0 725--750, 2014.

\bibitem[Sirki\"{a} et~al.(2007)Sirki\"{a}, Taskinen, and Oja]{sirkia07}
S.~Sirki\"{a}, S.~Taskinen, and H.~Oja.
\newblock Symmetrised {M}-estimators of scatter.
\newblock \emph{J. Mult. Anal.}, 98:\penalty0 1611--1629, 2007.

\bibitem[Soloveychik and Wiesel(2014)]{soloveychik14}
I.~Soloveychik and A.~Wiesel.
\newblock Performance {A}nalysis of {T}yler's {C}ovariance {E}stimator.
\newblock \emph{IEEE Trans. Sig. Process.}, 63-2:\penalty0 418--426, 2014.

\bibitem[Taskinen et~al.(2012)Taskinen, Koch, and Oja]{taskinen12}
S.~Taskinen, I.~Koch, and H.~Oja.
\newblock Robustifying principal component analysis with spatial sign vectors.
\newblock \emph{Statist. and Prob. Lett.}, 82:\penalty0 765--774, 2012.

\bibitem[Tukey(1975)]{tukey75}
J.W. Tukey.
\newblock Mathematics and picturing data.
\newblock In R.D. James, editor, \emph{Proceedings of the International
  Congress on Mathematics}, volume~2, pages 523--531, 1975.

\bibitem[Tyler(1987)]{tyler87}
D.E. Tyler.
\newblock A distribution-free {M}-estimator of multivariate scatter.
\newblock \emph{Ann. Statist.}, 15:\penalty0 234--251, 1987.

\bibitem[Visuri et~al.(2000)Visuri, Koivunen, and Oja]{visuri00}
S.~Visuri, V.~Koivunen, and H.~Oja.
\newblock Sign and rank covariance matrices.
\newblock \emph{J. Statist. Plan. Inf.}, 91:\penalty0 557--575, 2000.

\bibitem[Zuo(2003)]{zuo03}
Y.~Zuo.
\newblock Projection-based depth functions and associated medians.
\newblock \emph{Ann. Statist.}, 31:\penalty0 1460--1490, 2003.

\bibitem[Zuo and Cui(2005)]{ZuoCui05}
Y.~Zuo and M.~Cui.
\newblock Depth weighted scatter estimators.
\newblock \emph{Ann. Statist.}, 33-1:\penalty0 381--413, 2005.

\bibitem[Zuo and Serfling(2000)]{zuo00}
Y.~Zuo and R.~Serfling.
\newblock General notions of statistical depth functions.
\newblock \emph{Ann. Statist.}, 28-2:\penalty0 461--482, 2000.

\bibitem[Zuo et~al.(2004)Zuo, Cui, and He]{ZuoCuiHe04}
Y.~Zuo, M.~Cui, and X.~He.
\newblock On the {S}taehl-{D}onoho estimator and depth-weighted means of
  multivariate data.
\newblock \emph{Ann. Statist.}, 32-1:\penalty0 167--188, 2004.

\end{thebibliography}

\end{document}